\documentclass[UTF-8,reqno]{amsart}
\usepackage{enumerate}
\usepackage{verbatim,listings}
\usepackage{amssymb,url,color, booktabs}

\newcommand{\bee}{\begin{equation}}
\newcommand{\ee}{\end{equation}}
\newcommand{\bbb}{\begin{align}}
\newcommand{\eee}{\end{align}}

\usepackage{mathrsfs}
\usepackage{amsmath}
\usepackage{verbatim}
\usepackage{fancyhdr}
\usepackage[nobysame]{amsrefs}
\BibSpec{article}{%
+{}{\PrintAuthors} {author}
+{,}{ \textrm} {title}
+{.}{ \textit} {journal}
+{,}{ \textbf} {volume}
+{}{ \parenthesize} {date}
+{,}{ } {pages}
+{.}{ arXiv:} {eprint}
+{.}{} {transition}
}
\BibSpec{book}{%
+{}{\PrintAuthors} {author}
+{,}{ \textit} {title}
+{.}{ \textrm} {series} 
+{,}{ Vol.} {volume}
+{.}{ } {publisher}
+{,}{ } {date}
+{.}{} {transition}
}
\usepackage{color}
\usepackage{tikz}
\usepackage{graphicx}
\usepackage[colorlinks=true]{hyperref}
\hypersetup{
    linkcolor=blue,          
    citecolor=red,        
    filecolor=blue,      
    urlcolor=cyan
}

\newcommand{\red}{\color{red}}

\numberwithin{equation}{section}

\newtheorem{theorem}{Theorem}[section]

\newtheorem{lemma}[theorem]{Lemma}
\newtheorem{definition}[theorem]{Definition}
\newtheorem{proposition}[theorem]{Proposition}

\newtheorem{corollary}[theorem]{Corollary}

\newenvironment{proof of theorem 1.2 and 1.3}{{\it Proof of Theorem 1.2 and 1.3}.}{{\hfill 	
$\square$\hskip - \parfillskip}}
\newenvironment{proof of theorem 1.4}{{\it Proof of Theorem 1.4}.}{{\hfill 	
		$\square$\hskip - \parfillskip}}
\newenvironment{proof of theorem 1.5}{{\it Proof of Theorem 1.5}.}{{\hfill 	
		$\square$\hskip - \parfillskip}}
\newenvironment{proof of theorem 1.6}{{\it Proof of Theorem 1.6}.}{{\hfill 	
		$\square$\hskip - \parfillskip}}

\makeatletter

\newcommand{\Rmnum}[1]{\expandafter\@slowromancap\romannumeral #1@}
\makeatother

\usepackage{caption}
\def\geq{\geqslant}
\def\leq{\leqslant}
\def\ge{\geqslant}
\def\le{\leqslant}

\def\R{\mathbb{R}}
\def\Sn{\mathbb{S}^{n}}
\def\={&\!\!=\!\!&}
\def\al{\alpha}
\def\be{\beta}
\def\ga{\gamma}
\def\pa{\partial}
\def\lbd{\lambda}

\def\ep{\epsilon}

\def\Ctheta{\mathcal{C_{\theta}}}

\def\h{\hat{h}}
\def\p{\mathcal{P}}
\def\l{\mathcal{L}}
\def\r{\mathcal{R}}

\allowdisplaybreaks

	\tikzset{
	pattern size/.store in=\mcSize,
	pattern size = 5pt,
	pattern thickness/.store in=\mcThickness,
	pattern thickness = 0.3pt,
	pattern radius/.store in=\mcRadius,
	pattern radius = 1pt}
\makeatletter
\pgfutil@ifundefined{pgf@pattern@name@_ak5ydppfi}{
	\pgfdeclarepatternformonly[\mcThickness,\mcSize]{_ak5ydppfi}
	{\pgfqpoint{0pt}{0pt}}
	{\pgfpoint{\mcSize+\mcThickness}{\mcSize+\mcThickness}}
	{\pgfpoint{\mcSize}{\mcSize}}
	{
		\pgfsetcolor{\tikz@pattern@color}
		\pgfsetlinewidth{\mcThickness}
		\pgfpathmoveto{\pgfqpoint{0pt}{0pt}}
		\pgfpathlineto{\pgfpoint{\mcSize+\mcThickness}{\mcSize+\mcThickness}}
		\pgfusepath{stroke}
}}
\makeatother
\tikzset{every picture/.style={line width=0.75pt}} 

\begin{document}
\title{The capillary $L_{p}$ dual Christoffel-Minkowski-type problem for $p\geq q$ with $1\leq k\leq n$}\thanks{\it {The research is partially supported by NSFC (No. 12261105).}}

\author{Guanghan Li}
\address{School of Mathematics and Statistics\\                	
	Wuhan University\\                Wuhan 430072, P. R. China}
\email{ghli@whu.edu.cn}
\author{Mengliang Liu$^{*}$}\thanks{* Corresponding author}
\address{School of Mathematics and Statistics\\                	
	Wuhan University\\                Wuhan 430072, P. R. China}
\email{liumengliang@whu.edu.cn}

\begin{abstract}
We study the following fully nonlinear $k$-Hessian equation on a spherical cap $\mathcal{C}_{\theta}$ for $1\leq k\leq n$ and $0<\theta<\pi/2$:
$$
\left\{
\begin{aligned}
\sigma_k(\nabla^2h+h\sigma)
&=fh^{p-1}\bigl(h^2+|\nabla h|^2\bigr)^{\frac{k+1-q}{2}}
&&\text{in }\mathcal{C}_{\theta},\\
\nabla_\mu h&=\cot\theta\,h
&&\text{on }\partial\mathcal{C}_{\theta}.
\end{aligned}
\right.
$$
Under suitable convexity and boundary assumptions on $f$, we establish the existence of positive smooth strictly convex solutions for $p\geq q$, together with uniqueness when $p>q$ and $p\geq 1$. In the critical case $p=q>1$, scaling invariance prevents uniform $C^0$ estimates without normalization. We overcome this difficulty using a logarithmic gradient estimate and a normalized approximation with uniqueness up to dilation.
\end{abstract}

\maketitle
\setcounter{tocdepth}{2}

\section{Introduction}
The Brunn--Minkowski theory connects convex geometry with fully nonlinear elliptic equations. A basic inverse problem is the classical Minkowski problem. Its smooth form is a Monge--Amp\`ere equation on the unit sphere; see \cite{Minkowski1897,Aleksandrov1956,Nirenberg1953,Pogorelov1978,ChengYau1976,Schneider2014}. Christoffel--Minkowski problems form a related family. For the support function $h$ of a smooth strictly convex body in $\R^{n+1}$, their equation is
\begin{equation}\label{eq:classical-CM}
\sigma_k\bigl(\nabla^2h+h g^{\Sn}\bigr)=\varphi
\quad\text{on }\Sn,
\qquad 1\leq k\leq n.
\end{equation}
The cases $k=1$ and $k=n$ are the classical Christoffel and Minkowski problems \cite{Christoffel1865,Minkowski1897}. For $1<k<n$, Guan and Ma \cite{GuanMa2003} proved existence and strict convexity under a spherical convexity condition on $\varphi$. Related results can be found in \cite{GuanLinMa2006,GuanMaZhou2006,HuMaShen2004,GuanXia2018}.

Firey's $L_p$ addition and Lutwak's $L_p$ Brunn--Minkowski theory led to the $L_p$ Minkowski and $L_p$ Christoffel--Minkowski problems \cite{Firey1967,Firey1968,Lutwak1993}. The smooth equation is
\begin{equation}\label{eq:Lp-CM}
\sigma_k\bigl(\nabla^2h+h g^{\Sn}\bigr)=h^{p-1}\varphi
\quad\text{on }\Sn.
\end{equation}
Its solvability depends on both $p$ and $k$. We refer to \cite{GuanMa2003,GuanLinMa2006,GuanMaZhou2006,HuMaShen2004,g1,g2} and the references therein.

The dual Brunn--Minkowski theory gives another extension. Huang, Lutwak, Yang, and Zhang introduced dual curvature measures in \cite{HuangLutwakYangZhang2016}. Lutwak, Yang, and Zhang then introduced their $L_p$ version in \cite{LutwakYangZhang2018}. In the smooth closed setting, the related $L_p$ dual Christoffel--Minkowski equation is
\begin{equation}\label{eq:closed-Lp-dual-CM}
\sigma_k\bigl(\nabla^2h+h g^{\Sn}\bigr)
=fh^{p-1}\bigl(h^2+|\nabla h|^2\bigr)^{\frac{k+1-q}{2}}
\quad\text{on }\Sn.
\end{equation}
Here $p$ is the $L_p$ parameter and $q$ is the dual parameter. If $(p,q)=(1,k+1)$, Eq. \eqref{eq:closed-Lp-dual-CM} becomes \eqref{eq:classical-CM}. If $q=k+1$, it becomes \eqref{eq:Lp-CM}. If $k=n$, it is the smooth $L_p$ dual Minkowski equation.

The closed $L_p$ dual Minkowski problem has been studied in \cite{BoroczkyFodor2019,BoroczkyLutwakYangZhangZhao2017,ChenHuangZhao2019,ChenLi2021,HuangZhao2018,LiLiuLu2019,Zhao2017}; see also \cite{HuangLutwakYangZhang2016,LutwakYangZhang2018,matheng}. For general $k$, Ding and the first author \cite{DingLi2023} used inverse curvature flows to study Eq. \eqref{eq:closed-Lp-dual-CM}. Chen, Tu, and Xiang \cite{ChenTuXiang2025} treated $p\geq q$ with $p\geq1$. They proved existence and uniqueness when $p>q$. When $p=q>1$, they obtained an eigenvalue problem with uniqueness up to dilation. In the range $1<p<q\leq k+1$, Cabezas-Moreno and Hu \cite{CabezasMorenoHu2025} found even smooth strictly convex solutions by gradient estimates and degree theory. The top-order case $k=n$ is of Monge--Amp\`ere type. Its different parameter ranges lead to different existence and uniqueness results \cite{LutwakYangZhang2018,BoroczkyFodor2019,BoroczkyLutwakYangZhangZhao2017,ChenHuangZhao2019,ChenLi2021,HuangZhao2018,LiLiuLu2019,Zhao2017}.

For capillary hypersurfaces, the support function is defined on a spherical cap. Its natural boundary condition is of Robin type. Mei, Wang, and Weng introduced and solved the smooth capillary Minkowski problem \cite{MeiWangWeng2025Minkowski}. They also studied the capillary $L_p$ Minkowski problem \cite{MeiWangWeng2025LpMinkowski}. More generally, they defined the capillary $k$-th $L_p$ surface area measure and posed the capillary $L_p$ Christoffel--Minkowski problem \cite{prescrib}. Its smooth equation is
\begin{equation}\label{tt}
		\renewcommand{\arraystretch}{1.5}
		\left\{
		\begin{array}{rll}
			\sigma_k\big(\nabla^{2}h+h\sigma\big)&=fh^{p-1},& \text{in }\Ctheta,\\
		\nabla_{\mu}h&=\cot\theta\, h,& \text{on }\partial\Ctheta,
			\end{array}	
			\right.
	\end{equation}
where $\sigma$ is the standard metric on $\Ctheta$, $\mu$ is the outward unit conormal to $\partial\Ctheta$, and $\Ctheta$ is the spherical cap defined in Section \ref{sec 2}. A function $f:\Ctheta\to\mathbb{R}$ is called capillary even if
\begin{equation}\label{even}
f(\xi_1,\cdots,\xi_n,\xi_{n+1})=f(-\xi_1,\cdots,-\xi_n,\xi_{n+1})\quad\forall\xi\in\Ctheta.
\end{equation}
We say that a capillary hypersurface is capillary even if its support function is a capillary even function on $\Ctheta$. When $p=1$ and $1\leq k\leq n-1$, this is the capillary Christoffel--Minkowski equation. Hu, Ivaki, and Scheuer \cite{HuIvakiScheuer2025} and Mei, Wang, and Weng \cite{mwwcvpde} obtained existence results under different assumptions. The latter work assumes a suitable homotopy path. For $1<p<k+1$, related existence and uniqueness results for even data were proved in \cite{HuIvaki2026LpCurvature,prescrib}. Curvature flows have also been used in this setting \cite{HuHuIvaki2025Flow}.

Motivated by the analytic form of \eqref{eq:closed-Lp-dual-CM} and \eqref{tt}, we study the following Robin boundary value problem:
\begin{equation}\label{eq:capillary-Lp-dual-CM}
		\renewcommand{\arraystretch}{1.5}
		\left\{
		\begin{array}{rll}
			\sigma_k\big(\nabla^{2}h+h\sigma\big)&=fh^{p-1}\big(h^2+|\nabla h|^2\big)^{\frac{k+1-q}{2}},& \text{in }\Ctheta,\\
		\nabla_{\mu}h&=\cot\theta\, h,& \text{on }\partial\Ctheta.
			\end{array}	
			\right.
	\end{equation}
All derivatives and norms in \eqref{eq:capillary-Lp-dual-CM} are taken with respect to $\sigma$.  We emphasize that our formulation is at the level of the fully nonlinear PDE. The term ``dual Christoffel--Minkowski-type'' refers to the factor $\big(h^2+|\nabla h|^2\big)^{(k+1-q)/2}$ inherited from the closed Eq. \eqref{eq:closed-Lp-dual-CM}.

When $q=k+1$, the gradient factor disappears. Eq. \eqref{eq:capillary-Lp-dual-CM} then reduces to the capillary $L_p$ Christoffel--Minkowski equation studied in \cite{HuIvaki2026LpCurvature,prescrib}. For comparison, the formal case $k=n$ gives a Monge--Amp\`ere equation with a Robin boundary condition. Gao \cite{Gao2025CapillaryDual} proved existence and uniqueness for $p>q$ under the extra condition $q\leq1$. This condition was later removed in higher dimensions \cite{Gao2026HigherDimensions}. Hu, Hu, and Yang \cite{HuHuYang2026John} developed a capillary John ellipsoid theorem and a non-collapsing estimate. They established existence and uniqueness for $p>q$, obtained an eigenvalue result for $p=q>1$, and proved existence for $1<p<q<\infty$ using a degree-theoretic argument.
 For $1\leq k\leq n$, we studied in \cite{liliu} the Robin Hessian Eq. \eqref{eq:capillary-Lp-dual-CM} when $1<p<q\leq k+1$, and found an even smooth strictly convex solution under convexity and boundary conditions on $f$.

In the closed Christoffel--Minkowski theory, conditions such as
$$
\nabla^2\varphi^{-1/k}+g^{\Sn}\varphi^{-1/k}\geq0
$$
are naturally connected with constant-rank arguments. In the $L_p$-Christoffel--Minkowski problem, the corresponding exponent becomes $-(p+k-1)^{-1}$.  In the spherical cap setting, we combine this interior condition with the following Robin-type boundary condition:
$$
\nabla^2 f^{-\frac{1}{p+k-1}}
+\sigma f^{-\frac{1}{p+k-1}}\ge0
\quad\text{in }\Ctheta,
\qquad
\nabla_{\mu}f^{-\frac{1}{p+k-1}}
\le\cot\theta f^{-\frac{1}{p+k-1}}
\quad\text{on }\partial\Ctheta .
$$
Conditions of this form appear in the capillary Christoffel--Minkowski and capillary $L_p$ Christoffel--Minkowski problems \cite{HuIvaki2026LpCurvature,HuIvakiScheuer2025,mwwcvpde,prescrib}. These conditions ensure that the spherical Hessian remains of full rank when passing to the limit.

We say that a function $h\in C^{2}(\Ctheta)$ is convex (respectively, strictly convex) if
$$b[h]:=\nabla^{2}h+h\sigma$$
is positive semidefinite (respectively, positive definite). When $\theta=\frac{\pi}{2}$, a reflection argument reduces the problem to the classical closed case. We therefore focus on the general case $\theta\neq\frac{\pi}{2}$.

\begin{theorem}\label{thm1}
Let $1\leq k\leq n-1$ and $0<\theta<\pi/2$. Assume that $f\in C^{\infty}(\Ctheta)$ is positive and satisfies
\begin{equation}\label{23}
\nabla^{2}{f}^{-\frac{1}{p+k-1}}+\sigma f^{-\frac{1}{p+k-1}}\geq 0\quad\rm{in}\quad\Ctheta,
\end{equation}
and
\begin{equation}\label{24}
\nabla_{\mu}f^{-\frac{1}{p+k-1}}\leq \cot\theta f^{-\frac{1}{p+k-1}}\quad\rm{on}\quad\partial\Ctheta.
\end{equation}
\begin{enumerate}
\item If $p>q$ and $p\geq1$, then Eq. \eqref{eq:capillary-Lp-dual-CM} has a unique positive smooth strictly convex solution.
\item If $p=q>1$, then there is a unique positive constant $\tau$ such that
\begin{equation*}
		\renewcommand{\arraystretch}{1.5}
		\left\{
		\begin{array}{rll}
			\sigma_k\big(\nabla^{2}h+h\sigma\big)&=\tau fh^{p-1}\big(h^2+|\nabla h|^2\big)^{\frac{k+1-q}{2}},& \rm{in}\quad\Ctheta,\\
		\nabla_{\mu}h&=\cot\theta\, h,& \rm{on}\quad\partial\Ctheta,
			\end{array}	
			\right.
	\end{equation*}
has a positive smooth strictly convex solution. This solution is unique up to a positive dilation.
\end{enumerate}
\end{theorem}

For $k=n$, Hu, Hu, and Yang \cite{HuHuYang2026John} previously established part (1) of Theorem \ref{thm1} without any assumptions on $f$, while their result corresponding to part (2) requires $f$ to be capillary even.
In contrast, our approach establishes both conclusions for $k=n$ without imposing any such restriction on $f$.

\begin{theorem}\label{thm11}
Let $0<\theta<\pi/2$, and $f\in C^{\infty}(\Ctheta)$ be a positive smooth function. 
\begin{enumerate}
\item If $p>q$ and $p\geq1$, then Eq. \eqref{eq:capillary-Lp-dual-CM} has a unique positive smooth strictly convex solution for $k=n$.
\item If $p=q>1$, then there is a unique positive constant $\tau$ such that
\begin{equation*}
		\renewcommand{\arraystretch}{1.5}
		\left\{
		\begin{array}{rll}
			\sigma_n\big(\nabla^{2}h+h\sigma\big)&=\tau fh^{p-1}\big(h^2+|\nabla h|^2\big)^{\frac{n+1-q}{2}},& \rm{in}\quad\Ctheta,\\
		\nabla_{\mu}h&=\cot\theta\, h,& \rm{on}\quad\partial\Ctheta,
			\end{array}	
			\right.
	\end{equation*}
has a positive smooth strictly convex solution. This solution is unique up to a positive dilation.
\end{enumerate}
\end{theorem}

The main difficulty comes from the dual factor $\rho^{k+1-q}$, which does not appear in \eqref{tt}. Here
$$
\rho^2:=h^2+|\nabla h|^2.
$$
Differentiating it produces extra terms involving both $\nabla h$ and $\nabla^2h$. These terms enter the linearized equation. Near the boundary, they must be handled together with the Robin condition. Thus the closed-manifold estimates do not apply directly.

We first consider the case $p>q$. Although the maximum principle yields a $C^0$ estimate, it does not directly provide a gradient estimate, since an admissible solution is only $k$-convex at this stage and need not be convex. The cases $p>1$ and $p=1$ require different auxiliary functions. To handle the Robin boundary condition, we introduce a suitable boundary correction that allows us to apply the maximum principle.

The critical case $p=q>1$ has a different problem. The equation is invariant under dilation, so an unnormalized $C^0$ estimate is impossible. Following ideas from \cite{prescrib,MeiWangWeng2025LpMinkowski,d-g-l-l}, we prove the scale-invariant estimate
$$
\|\nabla\log h\|_{C^0(\Ctheta)}\leq C.
$$
The usual closed-manifold test function does not fit the Robin condition. We use a boundary correction inspired by the Neumann estimates in \cite{m-x,m-q} and the capillary argument in \cite{prescrib}. More precisely, we apply the maximum principle to
$$
|\nabla\omega|^2
=\bigl|\nabla(\log h-\cot\theta\,\tilde\ell)\bigr|^2.
$$
where $\tilde{\ell}$ is the positive function defined in the proof of Lemma \ref{hxy2} and \ref{hxy3}. By restricting $\theta \in (0,\frac{\pi}{2})$ and distinguishing the relevant ranges of $q$, we can effectively apply the maximum principle and obtain the desired logarithmic gradient estimate.

Another major challenge is to establish the global $C^2$ estimate. The Robin boundary condition provides only first-order information in the normal direction and therefore does not directly yield a bound for the double-normal derivative. Moreover, the gradient-dependent factor introduces additional terms when the equation is differentiated. To overcome these difficulties, we adapt the test function introduced in \cite{jfa} and incorporate the boundary barrier
$$
\mathcal{B}=\langle\nabla h,\nabla\zeta\rangle,
\qquad \zeta=1-e^{-d},
$$
where $d$ is a suitably chosen distance function on $\Ctheta$. This barrier allows us to control the boundary maximum, while the interior maximum principle yields the desired global $C^2$ estimate. Higher-order estimates then follow from standard elliptic regularity theory for uniformly elliptic equations with oblique boundary conditions. The full-rank theorem ensures that strict convexity is preserved when passing to the limit.
For $p>q$, the negative zeroth-order coefficient guarantees the invertibility of the linearized operator, thereby establishing openness in the continuity method. Finally, uniqueness follows from a comparison argument based on the quotient of two solutions.

The paper is organized as follows. Section \ref{sec 2} recalls basic facts about capillary convex hypersurfaces, the capillary Gauss map, capillary support functions, and elementary symmetric functions. Section \ref{sec 3} derives a priori estimates for solutions of \eqref{eq:capillary-Lp-dual-CM}, including gradient and logarithmic gradient estimates adapted to the Robin boundary condition. Higher-order estimates then follow from the regularity theory for fully nonlinear, uniformly elliptic second-order equations with oblique derivative boundary conditions. In Section \ref{sec 4}, we first establish a full-rank theorem for the fully nonlinear elliptic equation \eqref{eq:capillary-Lp-dual-CM} subject to the Robin boundary condition. Finally, combining the full-rank theorem with the continuity method and an approximation argument, we complete the proofs of our main results.
 Throughout the paper, let $\{e_i\}_{i=1}^n$ be an orthonormal frame on $\Ctheta$. We write $h_{ij}=\nabla^2h(e_i,e_j)$ and $b_{ij,k}=\nabla_{e_k}b_{ij}$. We use similar notation for higher derivatives. Repeated indices are summed unless stated otherwise. The constants $C$ and $c$ may change from line to line.

\section{Preliminaries}\label{sec 2}
In this section, we recall some basic facts about capillary convex bodies and the $k$-th elementary symmetric function.
\subsection{Basics of capillary geometry}\label{sec2.1}
Let $\{E_i\}_{i=1}^{n+1}$ be the standard orthonormal basis of $\R^{n+1}$, and let $\R^{n+1}_{+}=\{x\in\R^{n+1}|  \langle x,E_{n+1}\rangle>0\}$ denote the upper Euclidean half-space. Let $\Sigma\subset\overline{\R^{n+1}_{+}}$ be a properly embedded, smooth, compact hypersurface with boundary such that
	$$\mathrm{int}(\Sigma)\subset\R^{n+1}_{+}\qquad \mathrm{and} \qquad\partial\Sigma\subset\partial\R^{n+1}_{+}.$$
We call $\Sigma\subset\overline{\R^{n+1}_{+}}$ a {\emph{capillary hypersurface}} if it intersects $\partial\R^{n+1}_{+}$ at a constant contact angle $\theta\in (0,\pi)$. Let $\nu$ denote the outward unit normal (the Gauss map) of $\Sigma$ with respect to $\widehat{\Sigma}$, the bounded closed region in $\overline{\R^{n+1}_{+}}$ enclosed by $\Sigma$ and $\partial\R^{n+1}_{+}$. The contact angle $\theta$ is defined by
	$$\cos(\pi-\theta) = \langle \nu , \textbf{e}\rangle,\ \mathrm{ on }\ \Sigma\cap \partial\R^{n+1}_{+},$$
	where $\textbf{e}:= -E_{n+1}$ is the unit outward normal of $\partial\R^{n+1}_{+}$.  Let $\mu$ be the unit outward conormal of $\partial\Sigma$ in $\Sigma$ and let $\bar{\nu}$ be the unit normal of $\partial\Sigma\subset\partial{\mathbb{R}^{n+1}_{+}}$. Assume that \{$\nu,\mu$\} and \{$\bar{\nu},\textbf{e}$\} have the same orientation in the normal bundle of $\partial\Sigma\subset\overline{\mathbb{R}^{n+1}_{+}}$. Then the following relation holds:
\begin{equation*}\label{normal transform}
\left\{\begin{aligned}
    \nu=&-\cos\theta\, \textbf{e}+\sin\theta\, \bar{\nu},\\
	\mu=&\sin\theta \,\textbf{e}+\cos\theta \,\bar{\nu}.
	\end{aligned}\right.
\end{equation*}

	If $\Sigma$ is convex, then the Gauss image $\nu(\Sigma)$ of $\Sigma$ lies in the spherical cap
	$$\Sn_{\theta}:=\left\{x\in \Sn | \langle x, E_{n+1}\rangle \geq \cos\theta\right\}.$$
	Instead of the usual Gauss map $\nu$, it is more convenient to use the following map
	$$\tilde{\nu}:= T \circ \nu: \Sigma\to \Ctheta,$$
	where $\Ctheta$ is the spherical cap defined by
	$$\Ctheta:= \left\{\xi\in\overline{\R^{n+1}_{+}} |~ |\xi-\cos\theta\, \textbf{e}|=1\right\}.$$
Here, $T: \Sn_{\theta}\to \Ctheta$ denotes the vertical translation $T(z)=z+\cos\theta\, \textbf{e}$. The diffeomorphism $\tilde{\nu}$ is called the {\emph{capillary Gauss map}} of $\Sigma$. We may therefore reparameterize $\Sigma$ by the inverse of this map on $\Ctheta$; see \cite{Xia-arxiv}*{Lemma 2.2}.
Setting $X=\tilde{\nu}^{-1}$, we define the support function $h = h_{\Sigma} : \mathcal{C}_\theta \to \mathbb{R}$ of a strictly convex capillary hypersurface $\Sigma$ by
\begin{equation*}
h(\xi) := \langle X(\xi) , \nu\big(X(\xi)\big)\rangle = \langle\tilde{\nu}^{-1}(\xi)  ,\xi - \cos \theta\, \textbf{e}\rangle, \ \forall \, \xi \in \mathcal{C}_\theta.
\end{equation*}

When $\Sigma=\mathcal{C}_{\theta}$, we have $X(\xi)=\xi$ for all $\xi\in \mathcal{C}_\theta$, and the support function of $\mathcal{C}_{\theta}$ is given by
\begin{equation*}
 \ell(\xi):=\langle\xi,\xi - \cos \theta\, \textbf{e}\rangle=\sin^2\theta+\cos \theta \langle \xi, \textbf{e} \rangle.
\end{equation*}
A direct calculation shows that $\nabla_{ij}^{2}\ell+\ell\delta_{ij}=\delta_{ij}$.
In place of the usual support function $h$, it is convenient to introduce the {\emph{capillary support function}}
$$u(\xi):=\frac{h(\xi)}{\ell(\xi)}.$$
The following properties of the support function $h$, the capillary support function $u$, and the spherical-cap support function $\ell$ will be used below. We refer to \cite{MeiWangWeng2025Minkowski} for the proof.
\begin{lemma}\label{lem1}
Along the boundary $\pa\Ctheta$, we choose an orthonormal frame $\{e_{i}\}_{i=1}^{n}$ with $e_{n}=\mu$, where $\mu$ is the unit outer normal of $\pa\Ctheta\subset \Ctheta$. The following boundary conditions hold on $\partial \mathcal{C}_\theta$:
\begin{enumerate}
\item $\nabla_{\mu}h=\cot\theta\, h$;
\item $h_{\al n}=0,\quad\al=1,2,\cdots,n-1$;
\item$ \nabla_{\mu}u=0$,\quad $\nabla_{\mu}\ell=\cot\theta\, \ell$;
\item $u_{\al n}=-\cot\theta\, u_{\al},\quad\al=1,2,\cdots,n-1$.
\end{enumerate}
\end{lemma}
\subsection{$k$-th elementary symmetric function}
Let $Z=\{z_{ij}\}$ be an $n\times n$ symmetric matrix, and define
$$
\sigma_{k}(Z)=\sigma_{k}(\Lambda(Z))=\sum_{1\leq i_{1}<i_{2}\cdots <i_{k}\leq n}\Lambda_{i_{1}}\Lambda_{i_{2}}\cdots \Lambda_{i_{k}},
$$
where $\Lambda:=\Lambda(Z)=(\Lambda_{1},\cdots,\Lambda_{n})\in \mathbb{R}^{n}$ is the vector of eigenvalues of $Z$.
\begin{definition}
For $1\leq k\leq n$, define the cone $\Gamma_{k}\subset\mathbb{R}^{n}$ by
$$
\Gamma_{k}=\{\Lambda \in \mathbb{R}^{n}:\sigma_{i}(\Lambda)>0, \ \forall 1\leq i \leq k\}.
$$
\end{definition}

\begin{definition}\label{kconvex}
A function $h\in C^{2}(\Ctheta)$ satisfying \eqref{eq:capillary-Lp-dual-CM} is called a $k$-convex solution if
$$
b[h]=\nabla^{2}h(\xi)+h(\xi)\sigma\in \Gamma_{k}
$$
for all points $\xi\in \Ctheta$. If $b[h]\in \Gamma_{n}$, then $h$ is strictly (spherical) convex.
\end{definition}
A function $h$ is called admissible if it is $k$-convex in the sense of Definition \ref{kconvex}. We denote by $\sigma_{k}(\Lambda|i)$ the symmetric function with $\Lambda_{i}=0$. We list below some basic properties of the $k$-th elementary symmetric functions that will be used throughout the paper.
\begin{proposition}\label{basic}
Let $\Lambda=(\Lambda_{1},\cdots, \Lambda_{n})\in \mathbb{R}^{n}$ and $k=0,1,\cdots, n$. Then,
\begin{enumerate}
\item $\sigma_{k+1}(\Lambda)=\sigma_{k+1}(\Lambda|i)+\lambda_{i}\sigma_{k}(\Lambda|i), \quad \forall 1\leq i\leq n$;
\item $\sum_{i=1}^{n}\Lambda_{i}\sigma_{k}(\Lambda|i)=(k+1)\sigma_{k+1}(\Lambda)$;

\item $\sum_{i=1}^{n}\sigma_{k}(\Lambda|i)=(n-k)\sigma_{k}(\Lambda)$;

\item $\sum^{n}_{i=1}\Lambda^{2}_{i}\sigma_{k-1}(\Lambda|i)=\sigma_{1}(\Lambda)\sigma_{k}(\Lambda)-(k+1)\sigma_{k+1}(\Lambda)$.
\end{enumerate}
\begin{proposition}\label{basic1}
For $1\leq k\leq n$ and $\Lambda\in \Gamma_{k}$, denote by $F(\Lambda):=(\sigma_{k}(\Lambda))^{\frac{1}{k}}$. If $\Lambda_{1}\geq \Lambda_{2}\geq \cdots\geq\Lambda_{n}$, then
\begin{enumerate}
\item $\frac{\pa F(\Lambda)}{\pa \Lambda_{1}}\leq  \frac{\pa F(\Lambda)}{\pa \Lambda_{2}}\leq \cdots\leq \frac{\pa F(\Lambda)}{\pa \Lambda_{n}}$;
\item $\sum_{i=1}^{n} \frac{\pa F(\Lambda)}{\pa \Lambda_{i}}\geq \binom{n}{k}^{\frac{1}{k}}$.
\end{enumerate}
\end{proposition}
\end{proposition}

\begin{proposition}[Newton-Maclaurin inequality]\label{NM} For $\Lambda\in \Gamma_{k}$ and $1\leq l \leq k\leq n$, we have
$$
\left[\frac{\sigma_{k}(\Lambda)}{\binom{n}{k}}\right]^{\frac{1}{k}}\leq \left[\frac{\sigma_{l}(\Lambda)}{\binom{n}{l}}\right]^{\frac{1}{l}}.$$
\end{proposition}
The following proposition is well known; see, e.g., \cite{Gcpam}*{Lemma 2.1}.
\begin{proposition}\label{g}
Denote by $\rm{Sym(n)}$ the set of all $n\times n$ symmetric matrices. Let $F$ be a $C^{2}$ symmetric function defined in some open subset $\Psi\subset \rm{Sym(n)}$. For a diagonal matrix $Z \in \Psi$ with distinct eigenvalues, write $F(Z)=f(\kappa(Z))$, where $\kappa(Z)=(\kappa_{1},\cdots,\kappa_{n})$ are the eigenvalues of $Z$. Let $\ddot{F}(Y,Y)$ be the second order derivative of a $C^{2}$ symmetric $F$ in direction $Y \in\rm{Sym(n)}$, then
$$\ddot{F}(Y,Y)=\sum_{j,k=1}^{n}\ddot{f}^{jk}Y_{jj}Y_{kk}+2\sum_{j<k}\frac{\dot{f}^{j}-\dot{f}^{k}}{\lambda_{j}-\lambda_{k}}Y^{2}_{jk},$$
where the first order derivative of $F$ satisfies $\dot{F}^{ij}(Z)=\dot{f}^{i}\delta_{ij}$.
\end{proposition}

\section{A priori estimates}\label{sec 3}
In this section, we derive a priori estimates for solutions of Eq. \eqref{eq:capillary-Lp-dual-CM}. The main result is the following.
\begin{theorem}\label{hxy66}
Let $1\leq k\leq n$ and $\theta\in(0,\frac{\pi}{2})$. For every positive smooth function $f$ on $\Ctheta$, the following statements hold:
\begin{enumerate}
\item If $p>q$, $p\geq1$, and the support function $h$ is an admissible solution to Eq. \eqref{eq:capillary-Lp-dual-CM} satisfying $b[h]\geq 0$, then
\begin{equation*}
\underset{\Ctheta}{\min}h\geq c,
\end{equation*}
and for any $\al\in(0,1)$,
\begin{equation}\label{hxy521}
||h||_{C^{4,\al}(\Ctheta)}\leq C,
\end{equation}
where the positive constants $c, C$ depend only on $n,p,q,\theta$ and $f$.
\item If $1<q<p\leq q+1$, then for every positive admissible solution $h$ to Eq. \eqref{eq:capillary-Lp-dual-CM} satisfying $b[h]\geq 0$, the rescaled function
$$\tilde{h}:=\frac{h}{\underset{\Ctheta}{\min}h}$$
satisfies
\begin{equation}\label{hxyhxy}
||\tilde{h}||_{C^{4,\al}(\Ctheta)}\leq C
\end{equation}
where $C>0$ depends only on $n,\theta$ and $f$, and is independent of $p$.
\end{enumerate}
\end{theorem}
\subsection{$C^{0}$ estimate} In this subsection, we derive a $C^{0}$ estimate for admissible solutions of Eq. \eqref{eq:capillary-Lp-dual-CM} when $p>q$.

\begin{lemma}\label{lem3}
Let $p>q$ and $\theta\in (0,\pi)$, and suppose that $h$ is a positive admissible solution of Eq. \eqref{eq:capillary-Lp-dual-CM}. Then there exists a positive constant $C$, depending on $n,k,p,q,\theta$, and $f$, such that
\begin{equation}\label{C0}
\frac{1}{C}\leq h \leq C.
\end{equation}
\end{lemma}
\begin{proof}
In terms of the capillary support function $u$, \eqref{eq:capillary-Lp-dual-CM} becomes
\begin{equation*}\label{u}
		\renewcommand{\arraystretch}{1.5}
		\left\{
		\begin{array}{rll}
			\sigma_{k}(\ell\nabla^2 u+\nabla u\otimes \nabla \ell +\nabla \ell\otimes \nabla u+u\sigma) &= f (u\ell)^{p-1}\rho^{k+1-q}, &{\rm in} \ \mathcal{C}_{\theta}, \\
		\nabla_{\mu}u&=0,& \text{on }\partial\Ctheta.
			\end{array}	
			\right.
	\end{equation*}
A direct computation shows that
$$ \min\{1-\cos\theta,\sin^{2}\theta\}\leq\ell\leq \max\{1-\cos\theta,\sin^{2}\theta\}.$$
Moreover,
$$ \min\{(1-\cos\theta)^{2},\sin^{2}\theta\}\leq\ell^{2}+|\nabla\ell|^{2}\leq \max\{(1-\cos\theta)^{2},\sin^{2}\theta\}.$$
By the boundary condition, we suppose that $u$ attains its maximum in $\Ctheta\setminus\pa\Ctheta$ at $\xi_{0}$. Then, at $\xi_{0}$,
$$\nabla u=0,\quad \nabla^{2}u\leq 0.$$
Therefore, we obtain
$$C_{n}^{k}u^{k}\geq f(u\ell)^{p-1}u^{k+1-q}(\ell^{2}+|\nabla \ell|^{2})^{\frac{k+1-q}{2}},$$
which implies
\begin{equation}\label{j1}
(\underset{\Ctheta}{\max}\ u)^{p-q}\leq C \frac{1}{\underset{\Ctheta}{\min}\ f},
\end{equation}
where $C$ is a positive constant that may change from line to line, and depends only on $n,k,p,q,\theta$.
Similarly, we also obtain
\begin{equation*}
(\underset{\Ctheta}{\min}\ u)^{p-q}\geq C\frac{1}{\underset{\Ctheta}{\max}\ f}.
\end{equation*}
\end{proof}
\subsection{$C^{1}$ estimates} In this subsection, we derive $C^{1}$ estimates for admissible solutions of Eq. \eqref{eq:capillary-Lp-dual-CM}. The cases $p>1$ and $p=1$ require different arguments and, consequently, different restrictions on the contact angle. We also establish a logarithmic gradient estimate for positive admissible solutions, which is crucial in the case $p=q>1$.

To prove the main theorem, we first introduce some notation. We rewrite \eqref{eq:capillary-Lp-dual-CM} as
\begin{equation}\label{y1}
F(b)=(\sigma_{k}(b))^{\frac{1}{k}}=f^{\frac{1}{k}}h^{\frac{p-1}{k}}\rho^{\frac{k+1-q}{k}}=\Psi(\xi,h,\nabla h).
\end{equation}
Define $$\Psi(\xi,z,y):=f(\xi)^{1/k}z^{\ga}(z^{2}+|y|^{2})^{\be/2},$$
where $\ga:=\frac{p-1}{k}$, $\be=\frac{k+1-q}{k}$. Moreover, we denote
$$F^{ij}:=\frac{\pa F}{\pa b_{ij}},\ \ \ F^{ij,kl}:=\frac{\pa^{2}F}{\pa b_{ij}\pa b_{ki}}.$$
We first derive the $C^{1}$ estimate for solutions of Eq. \eqref{eq:capillary-Lp-dual-CM} for $p>q$ and $p>1$. The construction of the auxiliary function is inspired by the ideas of Mei-Wang-Weng \cite{prescrib} for studying the Robin boundary problems. Define
$$d^*(\xi):=\text{dist}(\xi,\pa\Ctheta).$$
The function $d^*(\xi)$ is smooth near $\pa\Ctheta$; we extend it to a smooth function on $\Ctheta$ and still denote it by $d^*$. The function $d^*$ satisfies
$$d^*|_{\pa\Ctheta}=0,\quad\text{and}\quad\nabla d^*|_{\pa\Ctheta}=-\mu.$$
Let
$$
N(\xi):=1+\cot\theta\, d^*.
$$
For convenience, we can assume that
$$
\frac{1}{2}\leq N(\xi)\leq 2.
$$
\begin{lemma}\label{hxy1}
Let $p>q$, $p>1$, and $\theta\in (0,\pi)$. If $h$ is a positive admissible solution of Eq. \eqref{eq:capillary-Lp-dual-CM}, then
\begin{equation}\label{y2}
|\nabla h|\leq C,
\end{equation}
where $C>0$ depends on $\theta,n,k,p,q,\underset{\Ctheta}{\min}f$ and $||f||_{C^{1}(\Ctheta)}$.
\end{lemma}
\begin{proof}
Consider the auxiliary function
$$\mathcal{P}:=\frac{1}{2}\log (1+|\nabla \hat{h}|^{2})-\eta\log (S-\hat{h}-d^*)-a_{0}d^*,$$
where $\h:=Nh$ and $0<\eta\ll1$ to be determined later.
We choose the constants
\begin{equation}\label{y3}S:=2+2||h||_{C^{0}(\Ctheta)}+||d^*||_{C^{0}(\Ctheta)},\quad a_{0}:=\min\{\cot\theta,0\}.\end{equation}
It follows from \eqref{y3} that
\begin{equation}\label{b1}1\leq S-\h-d^*\leq2S\quad\text{in}\quad \Ctheta.\end{equation}
Lemma \ref{lem1} and the definition of $d^*$ imply that
\begin{equation}\label{y4}\nabla_{\mu}\h=\cot\theta\, h\nabla_{\mu}d^*+N\nabla_{\mu}h=0\quad \text{on}\quad \pa\Ctheta.\end{equation}
Assume $\p$ attains its maximum value at some point, say $\xi_{0}\in \Ctheta$. In view of the position of $\xi_{0}$, we divide the proof into the following two cases. We may further assume that $|\nabla \h|(\xi_{0})\gg1$, otherwise we are done.\\

\noindent
{\bf Case 1.} $\xi_{0}\in\pa\Ctheta$. Choose a local orthonormal frame $\{e_{i}\}_{i=1}^{n}$ around $\xi_{0}$ such that $e_{n}=\mu$ and $\nabla \h=|\nabla \h|e_{1}$ at $\xi_{0}$. To simplify the computations, define
$$
Q:=1+|\nabla \h|^{2},\quad \text{and}\quad D:=S-\h-d^{*}.
$$
By the Gauss-Weingarten formula for $\pa\Ctheta\subset\Ctheta$ and \eqref{y4}
$$\h_{n1}=\nabla_{e_{1}}(\nabla_{e_{n}}\h)-\sigma(\nabla_{e_{1}}e_{n},\nabla\h)=-\cot\theta\, \h_{1}.$$
The maximum principle on the boundary implies
\begin{align*}
0\leq \nabla_{\mu}\p&=\frac{\sum_{s=1}^{n}\h_{s}\h_{sn}}{Q}-\eta\frac{-\h_{\mu}-d^*_{\mu}}{D}+a_{0}\\
&=\frac{\h_{1}\h_{1n}}{Q}-\frac{\eta}{D}+a_{0}\\
&=-\Big(\frac{\cot\theta\,\h_{1}^{2}}{Q}+\frac{\eta}{D}\Big)+a_{0}<0,
\end{align*}
which is a contradiction.

\noindent
{\bf Case 2.} $\xi_{0}\in \Ctheta\setminus\pa\Ctheta$. Choose a local orthonormal frame $\{e_{i}\}_{i=1}^{n}$ around $\xi_{0}$ such that $b[h]=\nabla^{2}h+h\sigma$ is diagonal at $\xi_{0}$, and let $\lambda_{i}$ denote the eigenvalues of $b[h]$. Define the elliptic operator
$$
\mathcal{L}:=F^{ij}\nabla_{ij}-\Psi_{y_{m}}\nabla_{m},
$$
where $\Psi_{y_{m}}:=\frac{\Psi}{\rho^{2}}h_{m}\be$.

By the maximum value condition, we obtain
\begin{equation}\label{y5}
0=\p_{i}=\frac{\sum_{s=1}^{n}\h_{s}\h_{si}}{Q}+\frac{\eta(\h_{i}+d^*_{i})}{D}-a_{0}d^*_{i},
\end{equation}
and
\begin{equation}\label{y6}
0\geq \mathcal{L}\p=\frac{F^{ij}\h_{si}\h_{sj}+\h_{s}\mathcal{L}\h_{s}}{Q}-2F^{ij}U_{i}U_{j}+\eta\frac{\mathcal{L}\h+\l d^*}{D}+\eta F^{ij}B_{i}B_{j}-a_{0}\mathcal{L}d^*,
\end{equation}
where $U_{i}:=\frac{\h_{s}\h_{si}}{Q}$ and $B_{i}:=\frac{\h_{i}+d^*_{i}}{D}$.
Differentiating \eqref{y1} along the direction $e_{s}$, we obtain
\begin{equation}\label{y66}
F^{ij}h_{ijs}=\Psi_{\xi_{s}}+\Psi_{z}h_{s}+\Psi_{y_{m}}h_{ms}-h_{s}\sum_{i=1}^{n}F^{ii}.
\end{equation}
On the spherical cap $\Ctheta$, the commutator formula holds:
\begin{equation}\label{y7}
h_{kij}=h_{ijk}+h_{k}\delta_{ij}-h_{j}\delta_{ki}.
\end{equation}
Combining \eqref{y66} and \eqref{y7}, we derive
\begin{equation*}
F^{ii}h_{sii}=\Psi_{\xi_{s}}+\Psi_{z}h_{s}+\Psi_{y_{m}}h_{ms}-F^{ss}h_{s},
\end{equation*}
which implies
\begin{equation}\label{y8}
\mathcal{L}h_{s}=\Psi_{\xi_{s}}+\Psi_{z}h_{s}-F^{ss}h_{s}.
\end{equation}
By a direct computation, we obtain
\begin{equation}\label{y9}
\mathcal{L}\h_{s}=N\l h_{s}+h_{s}\l N+N_{s}\l h+h\l N_{s}+2F^{ii}N_{i}h_{si}+2F^{ii}h_{i}N_{si}.
\end{equation}
Substituting \eqref{y8} into \eqref{y9} gives
\begin{align*}
\frac{\h_{s}\l \h_{s}}{Q}=&\frac{N\h_{s}\Psi_{\xi_{s}}}{Q}+\frac{N\Psi_{z}h_{s}\h_{s}}{Q}-\frac{NF^{ss}h_{s}\h_{s}}{Q}+\frac{h_{s}\h_{s}\l N}{Q}\\
+&\frac{N_{s}\h_{s}\l h}{Q}+\frac{2F^{ii}N_{i}\h_{s}h_{si}}{Q}+\frac{2F^{ii}h_{i}\h_{s}N_{si}}{Q}+\frac{h\h_{s}\l N_{s}}{Q}.
\end{align*}
On the one hand, for $|\nabla \h|$ sufficiently large, there exists a positive constant $C$ such that $\frac{1}{C}|\nabla \h|\leq |\nabla h|\leq C|\nabla \h|$. Consequently,
\begin{equation*}
|\Psi_{y}|\leq \be\Psi\frac{|\nabla \h|}{\rho^{2}}\leq C\frac{\Psi}{|\nabla \h|} \quad\text{and}\quad \h_{s}h_{s}\geq C|\nabla \h|^{2}.
\end{equation*}
When $p>1$, we obtain
\begin{equation*}
\Psi_{z}=\frac{\Psi}{h}(\ga+\be\frac{h^{2}}{\rho^{2}})\geq \frac{\min\{\ga,\ga+\be\}}{\underset{\Ctheta}{\max}h}\Psi:=c^*_{1}\Psi.
\end{equation*}
On the other hand, we have
$$
|\l N|+\sum_{s}|\l N_{s}|\leq C\left(\sum_{i=1}^{n}F^{ii}+|\Psi_{y}|\right),
$$
and for sufficiently small \(\epsilon\),
$$
\frac{2F^{ii}N_{i}h_{si}\h_{s}}{Q}\geq-\frac{\epsilon}{Q}F^{ii}\lambda_{i}^{2}-C\sum_{i=1}^{n}F^{ii}-C\sum_{i=1}^{n}F^{ii}\frac{1}{|\nabla \h|},
$$
where we have used
$$
2F^{ii}|\h_{i}N_{i}\lambda_{i}|\leq\epsilon F^{ii}\lambda_{i}^{2}+C(\epsilon)\sum_{i=1}^{n}F^{ii}|\nabla \h|^{2}.
$$
In particular, we can estimate the term \(\frac{N\Psi_{z}\h_{s}h_{s}}{Q}\) as follows:
$$
\frac{N\Psi_{z}\h_{s}h_{s}}{Q}\geq\frac{Nc_{1}^{*}\Psi C|\nabla \h|^{2}}{Q}\geq c_{2}^{*}\Psi.
$$
Combining the above estimates, we derive
\begin{equation}\label{y16}
\frac{\h_{s}\l \h_{s}}{Q}\geq c_{2}^*\Psi-\frac{\ep}{Q}F^{ii}\lambda_{i}^{2}-C\sum_{i=1}^{n}F^{ii}-C\frac{\Psi}{|\nabla \h|}.
\end{equation}
Observe that
\begin{align}\label{y11}
\l \h&=N\l h+h\l N+2F^{ii}N_{i}h_{i}\notag\\
&=N\Psi -Nh\sum_{i=1}^{n}F^{ii}-Nh_{m}\Psi_{y_{m}}+hF^{ii}N_{ii}-hN_{m}\Psi_{y_{m}}+2F^{ii}N_{i}h_{i}\notag\\
&=N\Psi(1-\be\frac{|\nabla h|^{2}}{\rho^{2}})-Nh\sum_{i=1}^{n}F^{ii}+hF^{ii}N_{ii}-hN_{m}\Psi_{y_{m}}+2F^{ii}N_{i}h_{i}\notag\\
&\geq -c_{3}^*\Psi-C(1+|\nabla \h|)\sum_{i=1}^{n}F^{ii}-C\frac{\Psi}{|\nabla \h|},
\end{align}
and
\begin{equation}\label{y14}
\l d^*\geq C\sum_{i=1}^{n}F^{ii}-C\frac{\Psi}{|\nabla \h|}.
\end{equation}
From \eqref{y5}, we have
\begin{align*}
-2F^{ii}U_{i}U_{i}+\eta F^{ii}B_{i}B_{i}&=-2F^{ii}(-\eta B_{i}+a_{0}d^*_{i})^{2}+\eta F^{ii}B_{i}^{2}\\
&=-2F^{ii}\Big(\eta^{2}B_{i}^{2}+a_{0}^{2}(d^*_{i})^{2}-2\eta a_{0}B_{i}d^*_{i}\Big)+\eta F^{ii}B_{i}^{2}\\
&=\eta(1-2\eta)F^{ii}B_{i}^{2}+4\eta a_{0}F^{ii}B_{i}d^*_{i}-2a_{0}^{2}F^{ii}(d_{i}^*)^{2}.
\end{align*}
By applying the Cauchy-Schwarz inequality and choosing $0<\eta\leq\frac{1}{4}$, we derive
$$F^{ii}4|\eta a_{0}B_{i}d^*_{i}|\leq \frac{\eta(1-2\eta)}{2}F^{ii}B_{i}^{2}+C(\eta,\theta)F^{ii}(d^*_{i})^{2}.$$
Therefore,
\begin{equation}\label{y12}
-2F^{ii}U_{i}U_{i}+\eta F^{ii}B_{i}B_{i}\geq -C\sum_{i=1}^{n}F^{ii}.
\end{equation}
For the term $\frac{F^{ij}\h_{si}\h_{sj}}{Q}$ in \eqref{y6}, we derive
\begin{align*}
F^{ij}\h_{si}\h_{sj}\geq& N^{2}F^{ii}\lambda_{i}^{2}-C|\nabla \h||\lambda_{i}|\sum_{i=1}^{n}F^{ii}-C(1+|\nabla \h|^{2})\sum_{i=1}^{n}F^{ii}\\
\geq&\frac{N^{2}}{2}F^{ii}\lambda_{i}^{2}-C(1+|\nabla \h|^{2})\sum_{i=1}^{n}F^{ii}.
\end{align*}
Hence, we obtain
\begin{equation}\label{y13}
\frac{F^{ij}\h_{si}\h_{sj}}{Q}\geq \frac{N^{2}}{2Q}F^{ii}\lambda_{i}^{2}-C\sum_{i=1}^{n}F^{ii}.
\end{equation}
Moreover, \eqref{y11} and \eqref{y14} yield
\begin{equation}\label{y15}
\eta\frac{\l d^*+\l \h}{D}\geq -\frac{\eta c_{3}^*}{2}\Psi-C(1+|\nabla \h|)\sum_{i=1}^{n}F^{ii}-C\frac{\Psi}{|\nabla \h|}.
\end{equation}
We may assume that $0<\eta\leq \min\left\{\frac{1}{4},\frac{c_{2}^{*}}{c_{3}^{*}}\right\}$, and take \(\epsilon=\frac{N^{2}}{4}\) in \eqref{y16}. Substituting \eqref{y13}, \eqref{y12}, \eqref{y15} and \eqref{y16} into \eqref{y6}, we obtain
\begin{equation}\label{y17}
0\geq \frac{N^{2}}{4Q}F^{ii}\lambda_{i}^{2}+\frac{c_{2}^*}{2}\Psi-C_{2}|\nabla \h|\sum_{i=1}^{n}F^{ii}-C_{1}\frac{\Psi}{|\nabla \h|}.
\end{equation}
Define
$$\mathcal{T}:=\left\{\ i\  \Big|1\leq i\leq n, |\h_{i}|\geq\frac{|\nabla \h|}{\sqrt{n}}\right\}.$$
Then $\mathcal{T}\neq\emptyset$. Fix $i_{0}\in \mathcal{T}$. By \eqref{y5},
\begin{align}\label{y19}
\h_{i_{0}i_{0}}&=-Q\Big(\frac{\eta}{D}(1+\frac{d^*_{i_{0}}}{\h_{i_{0}}})-a_{0}\frac{d^*_{i_{0}}}{\h_{i_{0}}}\Big)-\frac{1}{\h_{i_{0}}}\sum_{i\neq i_{0}}\h_{i}\h_{ii_{0}}\notag\\
&=-Q\Big(\frac{\eta}{D}(1+\frac{d^*_{i_{0}}}{\h_{i_{0}}})-a_{0}\frac{d^*_{i_{0}}}{\h_{i_{0}}}\Big)-\frac{1}{\h_{i_{0}}}\sum_{i\neq i_{0}}\h_{i}(N_{i}h_{i_{0}}+N_{i_{0}}h_{i}+hN_{ii_{0}}).
\end{align}
Observe that
\begin{align}\label{y20}
\Big|\frac{1}{\h_{i_{0}}}\sum_{i\neq i_{0}}\h_{i}(N_{i}h_{i_{0}}+N_{i_{0}}h_{i}+hN_{ii_{0}})\Big|\leq C_{3}(1+|\nabla \h|),
\end{align}
where $C_{3}$ depends only on $n$ and $||h||_{C^{0}(\Ctheta)}$. We may assume that
\begin{equation}\label{y21}
|\nabla \h|\geq \max\left\{6\sqrt{n},12S|a_{0}|\sqrt{n}/\eta,\frac{24S(2C_{3}+C_{4})}{\eta},2C_{1}/c_{2}^*,\frac{(2304nS^{2}+1)C_{2}}{\eta^{2}}\right\}.
\end{equation}
Together with \eqref{b1}, it yields
\begin{align}\label{y22}
\frac{\eta}{D}\Big(1+\frac{d^*_{i_{0}}}{\h_{i_{0}}}\Big)-\frac{|a_{0}|d^*_{i_{0}}}{\h_{i_{0}}}\geq \frac{\eta}{2S}\Big(1-\frac{\sqrt{n}|\nabla d^*|}{|\nabla \h|}\Big)-\frac{|a_{0}|\sqrt{n}|\nabla d^*|}{|\nabla \h|}\geq \frac{\eta}{3S}.
\end{align}
Substituting \eqref{y22} and \eqref{y20} into \eqref{y19}, and using \eqref{y21}, we infer that
\begin{align}\label{y23}
\lambda_{i_{0}}:=h_{i_{0}i_{0}}+h&\leq-\frac{\eta}{6S}Q+2C_{3}(1+|\nabla \hat{h}|)+C_{4}(1+|\nabla \hat{h}|)\notag\\
&\leq -\frac{\eta}{12S}Q.
\end{align}
Without loss of generality, assume that
$$\lambda_{1}\leq \lambda_{2}\leq \cdots\leq \lambda_{n}.$$
From \eqref{y23},
\begin{equation*}
\lambda_{1}\leq\lambda_{i_{0}}\leq-\frac{\eta}{12S}Q.
\end{equation*}
Substituting \eqref{y23} into \eqref{y17}, applying Proposition \ref{basic1}, we obtain
\begin{align*}
0\geq& \frac{1}{16Q}F^{11}\lambda_{1}^{2}-C_{2}|\nabla \h|\sum_{i=1}^{n}F^{ii}\\
\geq&\frac{\eta^{2}}{2304S^{2}n}\sum_{i=1}^{n}F^{ii}|\nabla \h|^{2}-C_{2}|\nabla \h|\sum_{i=1}^{n}F^{ii}>0.
\end{align*}
It follows that
$$\underset{\Ctheta}{\max}|\nabla \h|\leq C,$$
and therefore
$$\underset{\Ctheta}{\max}|\nabla h|\leq C,$$
where, in view of Lemma \ref{lem3}, the constant $C>0$ depends only on $p,q,n,k,\underset{\Ctheta}{\min}f$ and $||f||_{C^{1}(\Ctheta)}$. This completes the proof.
\end{proof}
As the proof of Lemma \ref{hxy1} shows, the restriction $p>1$ is required. We next consider the gradient estimate in the case $p=1>q$. The proof of the following lemma requires the angle to be restricted to $\theta \in(0,\frac{\pi}{2})$.
\begin{lemma}\label{hxy2}
Let $1=p>q$ and $\theta\in (0,\frac{\pi}{2})$. If $h$ is a positive admissible solution of Eq. \eqref{eq:capillary-Lp-dual-CM}, then
\begin{equation}\label{CC1}
|\nabla h|\leq C,
\end{equation}
where $C>0$ depends on $\theta,n,k,p,q,\underset{\Ctheta}{\min}f$ and $||f||_{C^{1}(\Ctheta)}$.
\end{lemma}
\begin{proof}
Define $v=\log h$. It follows from \eqref{eq:capillary-Lp-dual-CM} that $v$ satisfies
\begin{equation*}
		\renewcommand{\arraystretch}{1.5}
		\left\{
		\begin{array}{rll}
			\sigma_k(B)&=fe^{(1-q)v}(1+|\nabla v|^{2})^{\frac{k+1-q}{2}},& \text{in }\Ctheta,\\
		\nabla_{\mu}v&=\cot\theta ,& \text{on }\partial\Ctheta,
			\end{array}	
			\right.
	\end{equation*}
where $B_{ij}:=v_{ij}+v_{i}v_{j}+\sigma_{ij}$. We rewrite the equation as
\bee\label{xiang}
F(B)=\mathcal{R}(\xi,v,\nabla v)=f^{\frac{1}{k}}e^{\kappa v}(1+|\nabla v|^{2})^{\frac{\be}{2}},
\ee
where $\kappa:=\frac{1-q}{k}$, $\be=\frac{k+1-q}{k}=1+\kappa>1$.
Define
$$P:=|\nabla \omega|^{2},$$
where
$$\omega:=v-\cot\theta\,\tilde{\ell}=v-\phi,\quad \text{and}\quad\tilde{\ell}:=\frac{\ell}{\sin\theta\cos\theta}.$$
A direct computation shows that
$$\tilde{\ell}\Big|_{\pa\Ctheta}=\tan\theta,\quad \nabla\tilde{\ell}|_{\pa\Ctheta}=\mu,$$
and
\begin{equation}\label{lh1}
\nabla^{2}\tilde{\ell}=\frac{\nabla^{2}\ell}{\sin\theta\cos\theta}=\frac{1-\ell}{\sin\theta\cos\theta}\sigma,\quad\text{and}\quad \nabla^{2}\phi=\frac{1-\ell}{\sin^{2}\theta}\sigma.
\end{equation}
Suppose that $P$ attains its maximum at some point $\xi_{0}\in\Ctheta$. We again divide the proof into two cases: either $\xi_{0}\in \pa\Ctheta$ or $\xi_{0}\in \Ctheta\setminus\pa\Ctheta$.

\noindent
{\bf Case 1.} $\xi_{0}\in \pa\Ctheta$. Choose a local orthonormal frame $\{e_{i}\}_{i=1}^{n}$ around $\xi_{0}$ such that $e_{n}=\mu$. Using the properties of $\tilde{\ell}$, we have
\bee\label{grad-boundary-normal}
\omega_{n}=\nabla_{\mu}\omega=0\quad\text{on}\quad \pa\Ctheta.
\ee
By the Gauss-Weingarten equation for $\pa\Ctheta\subset\Ctheta$,
\bee\label{grad-boundary-mixed}
\omega_{\al n}=\nabla_{\al}(\omega_{n})-\sigma(\nabla\omega,\nabla_{e_{\alpha}}e_{n})=-\cot\theta\,\omega_{\al}.
\ee
Since $P$ attains its maximum at $\xi_{0}\in\pa\Ctheta$, it follows from \eqref{grad-boundary-normal} and \eqref{grad-boundary-mixed} that
\begin{align*}
0\leq\nabla_{n}P=2\sum_{i=1}^{n}\omega_{i}\omega_{in}=2\sum_{\al=1}^{n-1}\omega_{\al}\omega_{\al n}=-2\cot\theta\sum_{\al=1}^{n-1}\omega_{\al}^{2}.
\end{align*}
Since $\theta\in(0,\frac{\pi}{2})$, we obtain $|\nabla \omega|=0$. Therefore
\bee\label{ggb2}
|\nabla \log h|\leq \frac{1}{\sin^{2}\theta}\underset{\Ctheta}{\sup}|\nabla \ell|.
\ee
Combining this with \eqref{C0} proves \eqref{CC1}.

\noindent
{\bf{Case 2}.} $\xi_{0}\in\Ctheta\setminus\pa\Ctheta$. Choose an orthonormal frame $\{e_{i}\}_{i=1}^{n}$ around $\xi_{0}$, such that
$$B_{ij}:=v_{ij}+v_{i}v_{j}+\sigma_{ij}$$
is diagonal at $\xi_{0}$. Consequently,
$$F^{ij}:=\frac{\pa F(B)}{\pa B_{ij}}$$
is also diagonal at $\xi_{0}$. We may assume that $|\nabla v|$ is sufficiently large, which implies that $|\nabla \omega|$ is sufficiently large as well; otherwise, \eqref{CC1} follows immediately. Since $P$ attains its maximum at the interior point $\xi_{0}$, we have
\bee\label{lh4}
0=P_{i}=2\omega_{s}\omega_{si},
\ee
\bee\label{h5}
0\geq \frac{1}{2}F^{ij}P_{ij}=F^{ii}\omega_{si}^{2}+F^{ii}\omega_{s}\omega_{sii}.
\ee
For the standard metric on spherical cap $\Ctheta$, the commutator formula gives
\bee\label{h6}
v_{sij}=v_{ijs}+v_{s}\delta_{ij}-v_{j}\delta_{si}.
\ee
Differentiating Eq. \eqref{xiang} and using the fact $B_{ii,s}=v_{iis}+2v_{is}v_{i}$, we have
\bee\label{h7}
F^{ii}(v_{iis}+2v_{i}v_{si})=\mathcal{R}_{s}.
\ee
Combining \eqref{h5}, \eqref{h6} with \eqref{h7}, we derive
\begin{align}\label{h8}
0\geq\frac{1}{2}F^{ii}P_{ii}=&F^{ii}\omega_{si}^{2}+\omega_{s}\r_{s}-2F^{ii}v_{i}v_{is}\omega_{s}+\sum_{i=1}^{n}F^{ii}\omega_{s}v_{s}\notag\\
&-F^{ii}\omega_{i}v_{i}-F^{ii}\omega_{s}\phi_{sii}.
\end{align}
We next estimate the term $F^{ii}\omega_{si}^{2}-2F^{ii}v_{i}v_{is}\omega_{s}$ in \eqref{h8}. Observe that
\bee\label{h9}
\omega_{si}=v_{si}-\phi_{si}=B_{si}-v_{i}v_{s}-\delta_{si}-\phi_{si}=B_{si}-v_{i}v_{s}-\delta_{si} a_{1},
\ee
and
$$\omega_{s}v_{si}=\omega_{s}(\omega_{si}+\phi_{si})=a_{2}\omega_{i},$$
where we used \eqref{lh1} and \eqref{lh4} and $a_{1}:=1+\frac{1-\ell}{\sin^{2}\theta}$, $a_{2}:=\frac{1-\ell}{\sin^{2}\theta}$.
Therefore,
\begin{align}\label{h10}
F^{ii}\omega_{si}^{2}-2F^{ii}v_{i}v_{is}\omega_{s}&=F^{ii}\sum_{s}\Big((B_{si}-a_{1})\delta_{si}-v_{i}v_{s}\Big)^{2}-2a_{2}F^{ii}\omega_{i}v_{i}\notag\\
&=F^{ii}\Big(B_{ii}^{2}-2B_{ii}v_{i}^{2}+(|\nabla v|^{2}+2)v_{i}^{2}\Big)\notag\\
&\ \ \ \ \ \ \ -2a_{1}\r+a_{1}^{2}\sum_{i=1}^{n}F^{ii}+2a_{2}F^{ii}\phi_{i}v_{i},
\end{align}
where we used the identity $F^{ii}B_{ii}=\r$. Inserting \eqref{h10} into \eqref{h8} yields
 \begin{align}\label{h11}
 0\geq\frac{1}{2}F^{ii}P_{ii}&=F^{ii}\Big(B_{ii}^{2}-2B_{ii}v_{i}^{2}+(|\nabla v|^{2}+1)v_{i}^{2}\Big)+\sum_{i=1}^{n}F^{ii}|\nabla v|^{2}\notag\\
 &\ \ \ \ \ -2a_{1}\r+a_{1}^{2}\sum_{i=1}^{n}F^{ii}-\sum_{i=1}^{n}F^{ii} \langle \nabla \phi,\nabla v\rangle-F^{ii}\omega_{s}\phi_{sii}\notag\\
&\ \ \ \ \ +(2a_{2}+1)F^{ii}\phi_{i}v_{i}+\langle \nabla\omega,\nabla\r\rangle\notag\\
&\geq \sum_{i=1}^{n}F^{ii}|\nabla v|^{2}-2a_{1}\r-C|\nabla v|\sum_{i=1}^{n}F^{ii}+\langle \nabla\omega,\nabla\r\rangle,
 \end{align}
where we used that
$$B_{ii}^{2}-2B_{ii}v_{i}^{2}+(|\nabla v|^{2}+1)v_{i}^{2}=(B_{ii}-v_{i}^{2})^{2}+(|\nabla v|^{2}-v_{i}^{2})v_{i}^{2}+v_{i}^{2}\geq0.$$
From \eqref{xiang}, we obtain
\begin{align*}
\langle \nabla\omega,\nabla\r\rangle&=\r\Big(\langle \nabla(\log f^{1/k}),\nabla\omega\rangle+(\kappa+\frac{\be a_{2}}{1+|\nabla v|^{2}})\langle \nabla\omega,\nabla v\rangle\Big)\\
&\geq \r\Big(\frac{\kappa}{2}|\nabla v|^{2}-C_{5}|\nabla v|-C_{6}\Big).
\end{align*}
Since $\kappa>0$ is fixed, we may choose a constant $R_{1}=R_{1}(\kappa,C_{5},C_{6})$ sufficiently large such that, whenever $|\nabla v|\geq R_{1}$,
$$C_{5}|\nabla v|+C_{6}\leq \frac{\kappa}{4}|\nabla v|^{2}.$$
Therefore,
\bee\label{h12}
\langle \nabla\omega,\nabla\r\rangle\geq\frac{\kappa}{4}|\nabla v|^{2}\r.
\ee
Inserting \eqref{h12} into \eqref{h11}, we obtain
$$0\geq \sum_{i=1}^{n}F^{ii}|\nabla v|^{2}-2a_{1}\r+\frac{\kappa}{4}|\nabla v|^{2}\r-C_{7}|\nabla v|\sum_{i=1}^{n}F^{ii}.$$
Choose $R_{2}\geq R_{1}$ sufficiently large so that for $|\nabla v|\geq R_{2}$,
$$|\nabla v|^{2}-C_{7}|\nabla v|\geq \frac{|\nabla v|^{2}}{2},\quad \frac{\kappa}{4}|\nabla v|^{2}\r-2a_{1}\r\geq \frac{\kappa}{8}|\nabla v|^{2}\r.$$
This is a contradiction. Thus $|\nabla v|<R_{2}$ at the maximum point of $P$. Consequently,
$$|\nabla v|\leq C,$$
where the constant $C>0$ depends only on $p,q,n,k,\underset{\Ctheta}{\min}f$ and $||f||_{C^{1}(\Ctheta)}$. Together with the established $C^{0}$ estimate for $h$, this yields the desired global gradient estimate.
\end{proof}

We next establish a logarithmic gradient estimate for positive admissible solutions of Eq. \eqref{eq:capillary-Lp-dual-CM}.
\begin{lemma}\label{hxy3}
Let $p>q$, $p>1$, and $\theta\in (0,\frac{\pi}{2})$. If $h$ is a positive admissible solution of Eq. \eqref{eq:capillary-Lp-dual-CM}, then
\begin{equation}\label{logC1}
\underset{\Ctheta}{\max}~|\nabla \log h|\leq C,
\end{equation}
where $C>0$ depends on $\theta,n,k,p,q,\underset{\Ctheta}{\min}f$ and $||f||_{C^{1}(\Ctheta)}$. In particular, when $1<q< p\leq q+1$, the constant $C$ is independent of $p$.
\end{lemma}
\begin{proof}
Define $v=\log h$. It follows from \eqref{eq:capillary-Lp-dual-CM} that $v$ satisfies
\begin{equation*}
		\renewcommand{\arraystretch}{1.5}
		\left\{
		\begin{array}{rll}
			\sigma_k(B)&=fe^{(p-q)v}(1+|\nabla v|^{2})^{\frac{k+1-q}{2}},& \text{in }\Ctheta,\\
		\nabla_{\mu}v&=\cot\theta ,& \text{on }\partial\Ctheta,
			\end{array}	
			\right.
	\end{equation*}
where $B_{ij}:=v_{ij}+v_{i}v_{j}+\sigma_{ij}$. We rewrite the equation as
\begin{equation*}
F(B)=\tilde{\mathcal{R}}(\xi,v,\nabla v)=f^{\frac{1}{k}}e^{\tilde{\kappa} v}(1+|\nabla v|^{2})^{\frac{\be}{2}},
\end{equation*}
where $\tilde{\kappa}:=\frac{p-q}{k}\geq 0$, $\be:=\frac{k+1-q}{k}$. Consider the same auxiliary function as in Lemma \ref{hxy2}
$$P=|\nabla \omega|^{2},$$
where
$$\omega:=v-\cot\theta\,\tilde{\ell}=v-\phi,\quad \text{and}\quad\tilde{\ell}:=\frac{\ell}{\sin\theta\cos\theta}.$$

As in Lemma \ref{hxy2}, the auxiliary function $P$ attains its maximum at a point $\xi_{0}$ either on the boundary of the spherical cap or in its interior. In the boundary case, the gradient estimate follows directly from \eqref{ggb2}. In the interior case, the argument of Lemma \ref{hxy2} gives
\bee\label{x1}
0\geq |\nabla v|^{2}\sum_{i=1}^{n}F^{ii}+\langle \nabla \omega,\nabla \tilde{\r}\rangle -C|\nabla v|\sum_{i=1}^{n}F^{ii}-C\tilde{\r}.
\ee
By a direct computation, we have
\bee\label{x2}
\langle \nabla \omega,\nabla \tilde{\r}\rangle=\Big(\langle \nabla \log f^{1/k},\nabla \omega\rangle+(\tilde{\kappa}+\frac{\be a_{2}}{1+|\nabla v|^{2}})\langle \nabla v,\nabla \omega\rangle\Big)\tilde{\r}.
\ee
\noindent
{\bf Case 1.} $q\leq 1$. Then $\tilde{\kappa}=\frac{p-q}{k}\geq\frac{p-1}{k}:=\tilde{\kappa}_{1}$. We may choose a constant $H_{1}=H_{1}(\tilde{\kappa}_{1},C)$ sufficiently large such that, whenever $|\nabla v|\geq H_{1}$, we obtain
\bee\label{x3}
\langle \nabla \omega,\nabla \tilde{\r}\rangle\geq \frac{\tilde{\kappa}_{1}}{4}|\nabla v|^{2}\tilde{\r}.
\ee
Inserting \eqref{x3} into \eqref{x1}, we obtain
\bee\label{x4}
0\geq (|\nabla v|^{2}-C|\nabla v|)\sum_{i=1}^{n}F^{ii}+\frac{\tilde{\kappa}_{1}}{4}\tilde{\r}|\nabla v|^{2}-C\tilde{\r}.
\ee
Choose $H_{2}\geq H_{1}$ sufficiently large so that for $|\nabla v|\geq H_{2}$,
$$|\nabla v|^{2}-C|\nabla v|\geq \frac{|\nabla v|^{2}}{2},\quad \frac{\tilde{\kappa}_{1}}{4}|\nabla v|^{2}\tilde{\r}-C\tilde{\r}\geq \frac{\tilde{\kappa}_{1}}{8}|\nabla v|^{2}\tilde{\r}.$$
This is a contradiction in \eqref{x4} and thus completes the proof.

\noindent
{\bf Case 2.} $q>1$. Then $\tilde{\kappa}\geq 0$ and $\be=\frac{k+1-q}{k}<1$. We denote $\delta:=\max\{\be,0\}$; then we obtain
$$\tilde{\r}=f^{1/k}e^{\tilde{\kappa}v}(1+|\nabla v|^{2})^{\be/2}\leq C(1+|\nabla v|)^{\delta},$$
where $C>0$ depends only on $n,k,p,q,\underset{\Ctheta}{\min}\ f$ and $\underset{\Ctheta}{\max}\ f$. Moreover, for $q<p\leq q+1$, \eqref{j1} allows us to choose $C$ independently of $p$. From \eqref{x2}, we obtain
\bee\label{x5}
\langle \nabla \omega,\nabla \tilde{\r}\rangle\geq -C(1+|\nabla v|)\tilde{\r}\geq -C(1+|\nabla v|)^{1+\delta}.
\ee
In view of \eqref{x5} and \eqref{x1}, we obtain
\begin{align*}
0&\geq |\nabla v|^{2}\sum_{i=1}^{n}F^{ii}-C|\nabla v|\sum_{i=1}^{n}F^{ii}-C(1+|\nabla v|)^{1+\delta}\\
&\geq \frac{|\nabla v|^{2}}{2}\sum_{i=1}^{n}F^{ii}-C(1+|\nabla v|)^{1+\delta}\\
&\geq\binom{n}{k}^{\frac{1}{k}}\frac{|\nabla v|^{2}}{2}-C(1+|\nabla v|)^{1+\delta}>0.
\end{align*}
Consequently, we obtain the desired estimate \eqref{logC1}. Furthermore, if $1<q<p\leq q+1$, the constant $C$ in \eqref{logC1} is independent of $p$. This completes the proof of Lemma \ref{hxy3}.
\end{proof}
Motivated by \cite{Gcpam,jfa}, we consider the following general Christoffel--Minkowski-type equation with a Robin boundary condition:
\begin{equation}\label{general}
		\renewcommand{\arraystretch}{1.5}
		\left\{
		\begin{array}{rll}
			\sigma_{k}\big(\nabla^2h+h\sigma\big)&=\hat{f}(\xi, h, \nabla h),& \text{in }\Ctheta,\\
		\nabla_{\mu}h&=\cot\theta\, h,& \text{on }\partial\Ctheta.
			\end{array}	
			\right.
	\end{equation}
Let $\mathring{N}:=(0,\cdots,0,1-\cos\theta) \in \Ctheta$. Define $d:=\frac{1}{2\theta}d_{\mathring{N}}^{2}(\xi)$, where $d_{\mathring{N}}(\xi)$ denotes the geodesic distance from $\xi$ to $\mathring{N}$ on $\Ctheta$. Then $d$ is well defined and smooth on $\Ctheta$. A direct computation gives $d=\frac{\theta}{2}$ and $\nabla d=\mu$ on $\pa\Ctheta$, together with
$$(\nabla^{2}_{ij}d)\geq \min \left\{\frac{1}{\theta}, \cot\theta\right\}\delta_{ij},\quad \zeta_{\mu}|_{\pa\Ctheta}=e^{-\frac{\theta}{2}}.$$
Using this distance function, we construct the boundary barrier $\mathcal{B}:=\langle \nabla h, \nabla \zeta \rangle$, where $\zeta:=1-e^{-d}$. Suitable test functions and applications of the maximum principle on the boundary and in the interior then yield a global $C^{2}$ estimate.
\begin{theorem}\label{thm2}
Let $1\leq k\leq n-1$ and $\theta\in(0,\frac{\pi}{2})$. Suppose that $h$ is a positive, smooth, admissible solution of Eq. \eqref{general} satisfying $b[h]\geq 0$. Then there exists a constant $C>0$ such that
\bee\label{CCC2}\Delta h+nh\leq C,\ee
where $C$ depends on $n,k,p,q,\underset{\Ctheta}{\min}\ h, ||h||_{C^{1}(\Ctheta)}, \underset{\Ctheta}{\min} \hat{f}$ and $||\hat{f}||_{C^{2}(\Ctheta)}$. Furthermore, if $1<q< p\leq q+1$, the constant $C$ in \eqref{CCC2} is independent of $p$.
\end{theorem}
\begin{proof}
We consider the auxiliary function
\begin{equation*}
\Phi=\frac{1}{m}\log P_{m}+\frac{A}{2}\ell^{2}|\nabla u|^{2}-K\mathcal{B}-M\log h-L\zeta,
\end{equation*}
where $M,A,K,L>0$ are constants to be determined later, and
$$P_{m}= \sum_{j}\lambda_{j}^{m}, \quad m\geq 2.$$
Here $\lambda_{1}, \lambda_{2},\cdots, \lambda_{n}$ are the eigenvalues of the spherical Hessian $\nabla^{2}h+h\sigma$.
Assume that $\Phi$ attains its maximum at some point, say $\xi_{0} \in \Ctheta$. We divide the proof into two cases: either $\xi_{0} \in \pa\Ctheta$ or $\xi_{0} \in \Ctheta\setminus\pa\Ctheta$.

\noindent
{\bf Case I}. $\xi_{0} \in \pa\Ctheta$. Let $\{e_{i}\}_{i=1}^{n}$ be an orthonormal frame around $\xi_{0}$ such that $e_{n}=\mu$ and $b[h]$ is diagonal. First, we have
\begin{equation*}
\frac{1}{m}\nabla_{\mu}\log P_{m}=\frac{\sum_{\al}^{n-1}\lambda_{\al}^{m-1}b_{\al\al,\mu}+\lambda_{\mu}^{m-1}b_{\mu\mu,\mu}}{P_{m}}.
\end{equation*}
Differentiating \eqref{general} along the $\mu$-direction and using Lemma \ref{lem1} (2), we obtain
\begin{equation}\label{c22}
\hat{f}_{\mu}\leq C(1+b_{\mu\mu}).
\end{equation}
We also have
\begin{equation}\label{c23}
\sigma_{k}^{\mu\mu}b_{\mu\mu,\mu}+\sum_{\al}^{n-1}\sigma_{k}^{\al\al}b_{\al\al,\mu}=\hat{f}_{\mu}.
\end{equation}
Combining \eqref{c22} and \eqref{c23}, we obtain
\begin{align*}
\nabla_{\mu} \frac{1}{m}\log P_{m}&\leq \frac{\cot\theta(b_{\mu\mu}-b_{\al\al})\Big(\lambda_{\al}^{m-1}-\lambda_{\mu}^{m-1}\frac{\sigma_{k}^{\al\al}}{\sigma_{k}^{\mu\mu}}\Big)}{P_{m}}+\frac{\lambda_{\mu}^{m-1}\hat{f}C(1+b_{\mu\mu})}{P_{m}\sigma_{k}^{\mu\mu}}\\
&:=T_{1}+T_{2}.
\end{align*}
On the other hand, from Proposition \ref{basic} and \ref{NM}, we have
\begin{equation}\label{c24}
\hat{f}=\sigma_{k}=\lambda_{\mu}\sigma_{k-1}(\lambda|\mu)+\sigma_{k}(\lambda|\mu),
\end{equation}
and
\begin{equation}\label{c25}
\left[\frac{\sigma_{k-1}(\lambda|\mu)}{\binom{n-1}{k-1}}\right]^{\frac{1}{k-1}}\geq \left[\frac{\sigma_{k}(\lambda|\mu)}{\binom{n-1}{k}}\right]^{\frac{1}{k}}.
\end{equation}
Using \eqref{c25}, for some positive constant $C_{n,k}$, depending only on $n,k$, we get
$$\sigma_{k}(\lambda|\mu)\leq C_{n,k}\sigma_{k-1}(\lambda|\mu)\underset{i}{\max}\ \lambda_{i}.$$
So by \eqref{c24}, we have
$$\frac{\hat{f}}{\sigma_{k}^{\mu\mu}}\leq C\underset{i}{\max}\ \lambda_{i}.$$
Then
$$T_{2}:=\frac{\lambda_{\mu}^{m-1}\hat{f}C(1+b_{\mu\mu})}{P_{m}\sigma_{k}^{\mu\mu}}\leq\frac{C\underset{i}{\max}\ \lambda_{i}^{m}C(1+b_{\mu\mu})}{\underset{i}{\max}\ \lambda_{i}^{m}}:=C_{8}(1+b_{\mu\mu}).$$
At $\xi_{0}$, the maximum value condition implies
\begin{align*}
0\leq \nabla_{\mu}\Phi\leq& C(1+b_{\mu\mu})+A\ell\nabla_{\mu}\ell|\nabla u|^{2}+A\ell^{2}u_{k}u_{k\mu}-Ke^{-\frac{\theta}{2}}b_{\mu\mu}+KC-M\cot\theta-Le^{-\frac{\theta}{2}}\\
\leq & b_{\mu\mu}(C_{8}-Ke^{-\frac{\theta}{2}})+(C+KC-Le^{-\frac{\theta}{2}}),
\end{align*}
where we used Lemma \ref{lem1}. Choosing $L\gg K\gg1$ yields a contradiction.

\noindent
{\bf Case II}. $\xi_{0} \in \Ctheta\setminus\pa\Ctheta$. Hence at $\xi_{0}$, we have
\begin{equation}\label{c25-foc}
\frac{\sum_{j}\lambda_{j}^{m-1}b_{jj,i}}{P_{m}}+A\ell|\nabla u|^{2}\ell_{i}+A\ell^{2}u_{k}u_{ki}-K(h_{ki}\zeta_{k}+h_{k}\zeta_{ki})-M\frac{h_{i}}{h}-L\zeta_{i}=0,
\end{equation}
and
\begin{equation}
\begin{split}
\label{bjh}
0\geq &\frac{1}{P_{m}}\left(\sum_{j}\lambda^{m-1}_{j}b_{jj,ii}+(m-1)\sum_{j}\lambda^{m-2}_{j}b^{2}_{jj,i}+\sum_{p,q;p\neq q}\frac{\lambda^{m-1}_{p}-\lambda^{m-1}_{q}}{\lambda_{p}-\lambda_{q}}b^{2}_{pq,i}  \right)\\
&\quad -\frac{m}{P^{2}_{m}}(\sum_{j}\lambda^{m-1}_{j}b_{jj,i})^{2}-M\frac{h_{ii}}{h}+M\frac{h^{2}_{i}}{h^{2}}-L\zeta_{ii}+A\ell^{2}u_{ki}u_{ki}+A\ell^{2}u_{k}u_{kii}\\\
&+\ell_{i}^{2}A|\nabla u|^{2}+4A\ell u_{k}\ell_{i}u_{ki}-K(h_{kii}\zeta_{k}+2h_{ki}\zeta_{ki}+h_{k}\zeta_{kii})-A\ell^{2}|\nabla u|^{2}.\
\end{split}
\end{equation}
 Differentiating \eqref{general} twice, at $\xi_{0}$, we obtain
\begin{equation}\label{bjk}
\sum_{i}\sigma^{ii}_{k}b_{ii,j}=\hat{f}_{h_{j}}h_{jj}+\hat{f}_{h}h_{j}+\hat{f}_{j},
\end{equation}
and
\begin{equation}
\begin{split}
\label{bjl}
\sum_{i}\sigma^{ii}_{k}b_{ii,jj}+\sum_{p,q,r,s}\sigma^{pq,rs}_{k}b_{pq,j}b_{rs,j}&\geq  -C-C |h_{jj}|-C h^{2}_{jj}+\sum_{s}\hat{f}_{h_{s}}h_{sjj}\\
&\geq -C-Cb_{jj}-Cb^{2}_{jj}+\sum_{s}\hat{f}_{h_{s}}b_{sj,j}.
\end{split}
\end{equation}
Here, the positive constants $C$ depend on $\underset{\Ctheta}{\min}h, ||h||_{C^{1}(\Ctheta)}, \underset{\Ctheta}{\min}\hat{f}$ and $||\hat{f}||_{C^{2}(\Ctheta)}$.

Now using the Ricci identity $b_{ii,jj}=b_{jj,ii}+b_{ii}-b_{jj}$, we obtain
\begin{equation}\label{bjp}
\sum_{i}\sigma^{ii}_{k}b_{jj,ii}=\sum_{i}\sigma^{ii}_{k}b_{ii,jj}+b_{jj}\sum_{i}\sigma^{ii}_{k}-k\hat{f}.
\end{equation}
Multiplying both sides of \eqref{bjh} by $\sigma^{ii}_{k}$ and using \eqref{bjk}, \eqref{bjl} and \eqref{bjp}, we obtain

\begin{equation}
\begin{split}
\label{zx}
0&\geq \frac{1}{P_{m}}\sum_{j}\lambda^{m-1}_{j}\Big(-C_{7}-C_{8}b_{jj}-C_{9}b^{2}_{jj}+\sum_{s}\hat{f}_{h_{s}}b_{sj,j}-W(\sigma_{k})_{j}^{2}+W(\sigma_{k})_{j}^{2}\\
&+b_{jj}\sum_{i}\sigma^{ii}_{k}-k\hat{f}-\sum_{p,q,r,s}\sigma^{pq,rs}_{k}b_{pq,j}b_{rs,j}\Big)\\
&\quad +\frac{1}{P_{m}}(m-1)\sum_{i}\sigma^{ii}_{k}\sum_{j}\lambda^{m-2}_{j}b^{2}_{jj,i}+\frac{1}{P_{m}}\sum_{i}\sigma^{ii}_{k}\sum_{p,q;p\neq q}\frac{\lambda^{m-1}_{p}-\lambda^{m-1}_{q}}{\lambda_{p}-\lambda_{q}}b^{2}_{pq,i}\\
&\quad-\sum_{i}\frac{m\sigma^{ii}_{k}}{P^{2}_{m}}(\sum_{j}\lambda^{m-1}_{j}b_{jj,i})^{2}
-M\sum_{i}\sigma^{ii}_{k}\frac{h_{ii}}{h}+M\frac{\sum_{i}\sigma^{ii}_{k}h^{2}_{i}}{h^{2}}-L\sigma_{k}^{ii}\zeta_{ii}\\
&+A\sum_{i}\sigma_{k}^{ii}(\ell_{i}^{2}|\nabla u|^{2}+4\ell\ell_{i}u_{k}u_{ki}-\ell^{2}|\nabla u|^{2})+A\sum_{i}\sigma_{k}^{ii}(\ell^{2}u_{ki}u_{ki})+A\ell^{2}\sum_{i}\sigma_{k}^{ii}u_{k}u_{kii}\\
&-K\sum_{i}\sigma_{k}^{ii}(h_{kii}\zeta_{k}+2h_{ki}\zeta_{ki}+h_{k}\zeta_{kii}),
\end{split}
\end{equation}
where the positive constants $W$ are chosen sufficiently large, depending only on $\underset{\Ctheta}{\min}h$, $\|h\|_{C^1(\Ctheta)}$, $ \underset{\Ctheta}{\min}\hat{f}$, and $\|\hat{f}\|_{C^2(\Ctheta)}$, so that Lemma \ref{h1}-\ref{h3} holds.
We next estimate the terms $\frac{\sum_{s,j}\lbd_{j}^{m-1}\hat{f}_{h_ {s}}b_{sj,j}}{P_{m}}$ and $A\ell^{2}\sum_{i}\sigma_{k}^{ii}u_{k}u_{kii}$ in \eqref{zx}. Combining \eqref{c25-foc} with the identity
$$b_{ij}=\ell u_{ij}+\ell_{i}u_{j}+\ell_{j}u_{i}+u\delta_{ij},$$
we obtain
\begin{align}\label{c27}
&\frac{1}{P_{m}}\sum_{s,j}\lbd_{j}^{m-1}\hat{f}_{h_{s}}b_{sj,j}+A\ell^{2}\sum_{i}\sigma_{k}^{ii}u_{k}u_{kii}\notag\\
&=\sum_{s}\hat{f}_{h_{s}}\Big(-A\ell|\nabla u|^{2}\ell_{s}-A\ell^{2}u_{k}u_{ks}+K(h_{ks}\zeta_{k}+h_{k}\zeta_{ks})+M\frac{h_{s}}{h}+L\zeta_{s}\Big)\notag\\
&\quad\quad+A\sum_{i}\sigma_{k}^{ii}u_{k}\ell((\ell u_{ii})_{k}-\ell_{k}u_{ii})\notag\\
&= \sum_{s}\hat{f}_{h_{s}}\Big(-A\ell|\nabla u|^{2}\ell_{s}-A\ell^{2}u_{k}u_{ks}+K(h_{ks}\zeta_{k}+h_{k}\zeta_{ks})+M\frac{h_{s}}{h}+L\zeta_{s}\Big)\notag\\
&\quad\quad+A\sum_{i}\sigma_{k}^{ii}u_{k}\ell(b_{ii,k}-2\ell_{ik}u_{i}-2\ell_{i}u_{ik}-u_{k}-\ell_{k}u_{ii})\notag\\
&= \sum_{s}\hat{f}_{h_{s}}\Big(-A\ell|\nabla u|^{2}\ell_{s}-A\ell^{2}u_{k}u_{ks}+K(h_{ks}\zeta_{k}+h_{k}\zeta_{ks})+M\frac{h_{s}}{h}+L\zeta_{s}\Big)\notag\\
&\quad\quad+Au_{k}\ell\Big(\hat{f}_{k}+\hat{f}_{h}h_{k}+\hat{f}_{h_{s}}(\ell u_{sk}+u_{s}\ell_{k}+u_{k}\ell_{s}+u\delta_{sk}-\ell u\delta_{sk})\Big)-2A\ell \sum_{i}\sigma_{k}^{ii}u_{k}\ell_{ik}u_{i}\notag\\
&\quad\quad-2A\sum_{i}\sigma_{k}^{ii}u_{k}\ell_{i}(\ell u_{ik})-A\ell\sum_{i}\sigma_{k}^{ii}|\nabla u|^{2}-A\sum_{i}\sigma_{k}^{ii}u_{k}\ell_{k}(\ell u_{ii})\notag\\
&\quad\quad\geq -KC-AC-LC-MC-AC\sigma_{k}^{ii}-KC_{9}\underset{i}{\max}\lbd_{i}.
\end{align}
On the other hand,
\begin{equation}\label{zr}
-\sum_{p,q,r,s}\sigma^{pq,rs}_{k}b_{pq,j}b_{rs,j}=-\sum_{p,q}\sigma^{pp,qq}_{k}b_{pp,j}b_{qq,j}+\sum_{p,q}\sigma^{pp,qq}_{k}b^{2}_{pq,j}\ .
\end{equation}
Then substituting \eqref{c27} and \eqref{zr} into \eqref{zx} and using \eqref{bjk}, \eqref{zx} becomes
\begin{equation}
\begin{split}
\label{c26}
0&\geq \frac{1}{P_{m}}\sum_{j}\lambda^{m-1}_{j}\Big(-C_{0}(W)-C_{1}(W)b_{jj}-C_{2}(W)b^{2}_{jj}+W(\sigma_{k})_{j}^{2}\\
&-\sum_{p,q}\sigma^{pp,qq}_{k}b_{pp,j}b_{qq,j}+\sum_{p,q}\sigma^{pp,qq}_{k}b^{2}_{pq,j}\Big)+\frac{1}{P_{m}}(m-1)\sum_{i}\sigma^{ii}_{k}\sum_{j}\lambda^{m-2}_{j}b^{2}_{jj,i}\\
&\quad+\frac{1}{P_{m}}\sum_{i}\sigma^{ii}_{k}\sum_{p,q;p\neq q}\frac{\lambda^{m-1}_{p}-\lambda^{m-1}_{q}}{\lambda_{p}-\lambda_{q}}b^{2}_{pq,i}-\sum_{i}\frac{m\sigma^{ii}_{k}}{P^{2}_{m}}(\sum_{j}\lambda^{m-1}_{j}b_{jj,i})^{2}\\
&\quad +A\sum_{i}\sigma_{k}^{ii}(\ell^{2}u_{ki}u_{ki})+(-CK-CL-CA+M+M\frac{h_{i}^{2}}{h})\sum_{i}\sigma_{k}^{ii}\\
&\quad-\tilde{C}-C_{9}K\underset{i}{\max}\lbd_{i},
\end{split}
\end{equation}
where the positive constant $\tilde{C}$ depends on $\underset{\Ctheta}{\min}h, ||h||_{C^{1}(\Ctheta)}, \underset{\Ctheta}{\min}\hat{f}, ||\hat{f}||_{C^{2}(\Ctheta)}$ and constants $M,A,K,L$. Moreover, the constants $C_{0}(W), C_{1}(W), C_{2}(W)$ depend on $\underset{\Ctheta}{\min}h, ||h||_{C^{1}(\Ctheta)}, \underset{\Ctheta}{\min}\hat{f}, ||\hat{f}||_{C^{2}(\Ctheta)}$ and the positive constant $W$.

For the term $A\sum_{i}\sigma_{k}^{ii}(\ell^{2}u_{ki}u_{ki})$ in \eqref{c26}, we derive
\begin{align*}
A\sum_{i}\sigma_{k}^{ii}(\ell^{2}u_{ki}u_{ki})&=\sum_{i}\sigma_{k}^{ii}A(b_{ii}-2u_{i}\ell_{i}-u)^{2}\\
&\geq A\sum_{i}\sigma_{k}^{ii}(b_{ii}^{2}-4b_{ii}\ell_{i}u_{i}+4uu_{i}\ell_{i}-2ub_{ii})\\
&\geq A\sum_{i}\sigma_{k}^{ii}b_{ii}^{2}-CA-CA\sum_{i}\sigma_{k}^{ii}.
\end{align*}
Now, if we take $M\gg A\gg L\gg K\gg1$ while keeping the sign of the constant term $\tilde{C}$ unchanged, \eqref{c26} becomes
\begin{equation}
\begin{split}
\label{c27-final}
0&\geq \frac{1}{P_{m}}\sum_{j}\lambda^{m-1}_{j}(-C_{0}(W)-C_{1}(W)b_{jj}-C_{2}(W)b^{2}_{jj}+W(\sigma_{k})^{2}_{j}\\
&-\sum_{p,q}\sigma^{pp,qq}_{k}b_{pp,j}b_{qq,j}+\sum_{p,q}\sigma^{pp,qq}_{k}b^{2}_{pq,j})\\
&\quad+\frac{1}{P_{m}}(m-1)\sum_{i}\sigma^{ii}_{k}\sum_{j}\lambda^{m-2}_{j}b^{2}_{jj,i}+\frac{1}{P_{m}}\sum_{i}\sigma^{ii}_{k}\sum_{p,q;p\neq q}\frac{\lambda^{m-1}_{p}-\lambda^{m-1}_{q}}{\lambda_{p}-\lambda_{q}}b^{2}_{pq,i}\\
&\quad-\sum_{i}\frac{m\sigma^{ii}_{k}}{P^{2}_{m}}(\sum_{j}\lambda^{m-1}_{j}b_{jj,i})^{2}+A\sum_{i}\sigma^{ii}_{k}b^{2}_{ii}-\tilde{C}-C_{9}K\lbd_{1},
\end{split}
\end{equation}
where $\lbd_{1}:=\underset{i}{\max}\lbd_{i}$.
Next we deal with the third-order derivatives. Denote
$$\tilde{A}_{i}=\frac{\lambda^{m-1}_{i}}{P_{m}}(W(\sigma_{k})^{2}_{i}-\sum_{p,q}\sigma^{pp,qq}_{k}b_{pp,i}b_{qq,i}),$$
$$\tilde{B}_{i}=\frac{2}{P_{m}}\sum_{j}\lambda^{m-1}_{j}\sigma^{jj,ii}_{k}b^{2}_{jj,i},$$
$$\tilde{C}_{i}=\frac{m-1}{P_{m}}\sigma^{ii}_{k}\sum_{j}\lambda^{m-2}_{j}b^{2}_{jj,i},$$
$$\tilde{D}_{i}=\frac{2}{P_{m}}\sum_{j\neq i}\sigma^{jj}_{k}\frac{\lambda^{m-1}_{j}-\lambda^{m-1}_{i}}{\lambda_{j}-\lambda_{i}}b^{2}_{jj,i},$$
$$\tilde{E}_{i}=\frac{m\sigma^{ii}_{k}}{P^{2}_{m}}(\sum_{j}\lambda^{m-1}_{j}b_{jj,i})^{2}.$$

We now apply the key estimates from \cite{jfa}*{Lemmas 8 and 9}; see also \cite{CabezasMorenoHu2025}*{Lemmas 3.6 and 3.7}. The third-order derivative terms for $i=1$ and $i\neq1$ are treated using the strategy developed in \cite{Gcpam}.
\begin{lemma}\label{h1}
For any $i\neq 1$, we obtain
\[
\tilde{A}_{i}+\tilde{B}_{i}+\tilde{C}_{i}+\tilde{D}_{i}-\left( 1+\frac{1}{m} \right)\tilde{E}_{i}\geq 0,
\]
for sufficiently large $m$.
\end{lemma}

\begin{lemma}\label{h2}
For $l=1,\cdots,k-1$, if there exists a positive constant $\delta\leq 1$ such that $\lambda_{l}/\lambda_{1}\geq \delta$, then there exist two sufficiently small positive constants $\eta,\delta^{'}$ depending on $\delta$, such that, if $\lambda_{l+1}/\lambda_{1}\leq \delta^{'}$, we obtain
$$
\tilde{A}_{1}+\tilde{B}_{1}+\tilde{C}_{1}+\tilde{D}_{1}-\left(1+\frac{\eta}{m} \right)\tilde{E}_{1}\geq 0,
$$
for sufficiently large $m$.
\end{lemma}
\end{proof}
Arguing as in \cite[Corollary 10]{jfa}, we obtain the following result.
\begin{lemma}\label{h3}
There exist two finite sequences of positive numbers $\{\delta_{j}\}^{k}_{j=1}$ and $\{\xi_{j}\}^{k}_{j=1}$ such that, if the following inequality holds for some index $1\leq s\leq k-1$,
$$
\frac{\lambda_{s}}{\lambda_{1}}\geq \delta_{s}, \ {\rm and} \ \frac{\lambda_{s+1}}{\lambda_{1}}\leq \delta_{s+1},
$$
then for large $W$ in $\tilde{A}_{i}$, there holds
\begin{equation}\label{h4}
\tilde{A}_{1}+\tilde{B}_{1}+\tilde{C}_{1}+\tilde{D}_{1}-\left(1+\frac{\xi_{s}}{m} \right)\tilde{E}_{1}\geq 0.
\end{equation}
\end{lemma}
\begin{proof}
We construct the sequences $\{\delta_{j}\}_{j=1}^{k}$ and $\{\xi_{j}\}_{j=1}^{k}$ inductively. Set $\delta_{1} = 1/2$, so that $\lambda_{1}/\lambda_{1} = 1 > \delta_{1}$. The assertion for $j=1$ then follows directly from Lemma \ref{h2}.

Suppose that $\delta_{s}$ has been defined for $1 \leq s \leq k-1$. To determine $\delta_{s+1}$, apply Lemma \ref{h2} with $l=s$ and $\delta=\delta_{s}$. This gives a constant $\delta'_{s+1}>0$ such that, whenever $\lambda_{s+1} \leq \delta'_{s+1} \lambda_{1}$, estimate \eqref{h4} holds for some $\xi_{s}$. Choosing $\delta_{s+1} := \min\{\delta_{1}, \delta'_{s+1}\}$ ensures that \eqref{h4} holds whenever $\lambda_{s+1} \leq \delta_{s+1} \lambda_{1}$. This inductively determines both $\delta_{s+1}$ and $\xi_{s}$.
\end{proof}
\begin{proof}[\bf Continuation of the proof of Theorem \ref{thm2}.]
We distinguish two cases.

\noindent
{\bf Case I:} Suppose there exists an index $1 \leq s \leq k-1$ and positive constants $\{\delta_{j}\}_{j=1}^{k}$ such that
$$
\lambda_{s} \geq \delta_{s} \lambda_{1} \quad \text{and} \quad \lambda_{s+1} \leq \delta_{s+1} \lambda_{1}.
$$
Then, by combining Lemma \ref{h1} with Lemma \ref{h3}, we obtain
\begin{equation}\label{Pi}
\sum_{i}(\tilde{A}_{i}+\tilde{B}_{i}+\tilde{C}_{i}+\tilde{D}_{i})-\tilde{E}_{1}-\left(1+\frac{1}{m} \right)\sum^{n}_{i=2}\tilde{E}_{i}\geq 0.
\end{equation}

Recalling the definitions of $\tilde{A}_{i}$, $\tilde{B}_{i}$, $\tilde{C}_{i}$, $\tilde{D}_{i}$, and $\tilde{E}_{i}$ and substituting \eqref{Pi} into \eqref{c27-final}, we obtain
\begin{equation}
\begin{split}
\label{Qi}
0 &\geq \frac{1}{P_{m}}\sum_{j}\lambda^{m-1}_{j}\left(-C_{0}(W)-C_{1}(W)b_{jj}-C_{2}(W)b^{2}_{jj}\right) \\
&\quad + \sum^{n}_{i=2}\frac{\sigma^{ii}_{k}}{P^{2}_{m}}\left(\sum_{j}\lambda^{m-1}_{j}b_{jj,i}\right)^{2} + A\sum_{i}\sigma^{ii}_{k}b^{2}_{ii} - \tilde{C}-C_{19}K\lbd_{1},\\
&\geq -\frac{\tilde{C}_{0}(W)}{\lambda_{1}} - \tilde{C}_{1}(W) - \tilde{C} - \tilde{C}_{2}(W)\lambda_{1} + A\sigma^{11}_{k}b^{2}_{11}-C_{9}K\lbd_{1}.
\end{split}
\end{equation}
Applying the Newton--MacLaurin inequality once again yields
\begin{equation}\label{Mac}
\left[\frac{\sigma_{k-1}(\lambda|1)}{\binom{n-1}{k-1}}\right]^{\frac{1}{k-1}}
\geq
\left[\frac{\sigma_{k}(\lambda|1)}{\binom{n-1}{k}}\right]^{\frac{1}{k}}.
\end{equation}
Consequently, from \eqref{Mac}, there exists a positive constant $C_{n,k}$, depending only on $n$ and $k$, such that
\begin{equation}
\label{c28}
\sigma_{k}(\lambda|1)
\leq C_{n,k} \sigma_{k-1}(\lambda|1)^{\frac{k}{k-1}}
\leq C_{n,k} \lambda_{1} \sigma_{k-1}(\lambda|1),
\end{equation}
where the last inequality follows from $\lambda_1 \geq \lambda_j$ for all $j$.

Substituting \eqref{c28} into the decomposition $\sigma_{k}(\lambda)=\sigma_{k}(\lambda|1)+\lambda_{1}\sigma_{k-1}(\lambda|1)$, we obtain
\begin{equation*}
\lambda_{1}\sigma_{k-1}(\lambda|1) \geq C_{n,k} \sigma_{k}(\lambda).
\end{equation*}

Thus, we have
\begin{equation*}
\label{maxw11}
\frac{\sigma_{k}^{11} b_{11}^{2}}{\sigma_{k}}
= \frac{\sigma_{k-1}(\lambda|1) \lambda_{1}^{2}}{\sigma_{k}}
\geq \frac{C_{n,k} \sigma_{k} \lambda_{1}}{\sigma_{k}}
= C_{n,k} b_{11}.
\end{equation*}
Hence,
\begin{equation}\label{Ti}
\sigma^{11}_{k} b^{2}_{11} \geq C^{*} b_{11},
\end{equation}
where the positive constant $C^{*}$ depends on $n$, $k$, $\underset{\Ctheta}{\min}h$, $\|h\|_{C^{1}(\mathcal{C}_\theta)}$, $\underset{\Ctheta}{\min}\hat{f}$, and $\|\hat{f}\|_{C^{2}(\mathcal{C}_\theta)}$.

Inserting \eqref{Ti} into \eqref{Qi} and choosing
\[
A := \frac{\tilde{C}_{2}(W) + C_{9}K + 1}{C^{*}},
\]
we arrive at
\begin{equation}
\label{Ri}
\begin{aligned}
0 &\geq -\frac{\tilde{C}_{0}(W)}{\lambda_{1}}
      -\tilde{C}_{1}(W)
      -\tilde{C}
      -\tilde{C}_{2}(W)\lambda_{1}
      - C_{9}K \lambda_{1}
      + C^{*} A \lambda_{1} \\
&\geq -\frac{\tilde{C}_{0}(W)}{\lambda_{1}}
      -\tilde{C}_{1}(W)
      -\tilde{C}
      + \lambda_{1}.
\end{aligned}
\end{equation}
Thus, for sufficiently large $\lambda_{1}$, \eqref{Ri} implies
$$
\lambda_{1} \leq C
$$
for some positive constant $C$.

\noindent
{\bf Case II:} $\lambda_{k}\geq \delta_{k}\lambda_{1}$ with $\delta_{k}>0$. Since $\lambda_{1}\geq \lambda_{2}\geq \cdots \geq \lambda_{k}\geq \delta_{k}\lambda_{1}$,
$$
\hat{f}=\sigma_{k}\geq\lambda_{1}\cdots \lambda_{k}\geq \delta^{k-1}_{k}\lambda^{k}_{1},
$$
which again yields $\lambda_{1}\leq C$ for some constant $C>0$. Moreover, if $1<q<p\leq q+1$, the constant $C$ in \eqref{CCC2} is independent of $p$. This completes the proof of Theorem \ref{thm2}.
\end{proof}
Next, we establish the global $C^2$ estimate for solutions of \eqref{eq:capillary-Lp-dual-CM} in the case $k=n$. To this end, it suffices to show that Lemmas \ref{h1} and \ref{h2} also hold in this case.
\begin{corollary}\label{coro1}
Let $k=n$ and $\theta\in(0,\frac{\pi}{2})$. Suppose that $h$ is a positive, smooth, strictly convex solution of Eq. \eqref{general}. Then there exists a constant $C>0$ such that
\bee\label{CC2}\Delta h+nh\leq C,\ee
where $C$ depends on $n,k,p,q,\underset{\Ctheta}{\min}\ h, ||h||_{C^{1}(\Ctheta)}, \underset{\Ctheta}{\min} \hat{f}$ and $||\hat{f}||_{C^{2}(\Ctheta)}$. Furthermore, if $1<q<p\leq q+1$, the constant $C$ in \eqref{CC2} is independent of $p$.
\end{corollary}
\begin{proof}
We may assume that $\lbd_{1}\geq \lbd_{2}\geq \cdots\geq \lbd_{n}>0$. Recall that
$$\hat{f}:=\prod_{j=1}^{n}\lbd_{j}, \quad P_{m}:=\sum_{j}\lbd_{j}^{m}.$$
We also have that
$$\sigma_{n}^{jj}=\frac{\hat{f}}{\lbd_{j}},\quad
\sigma_{n}^{jj,ii}=
\begin{cases}
\frac{\hat{f}}{\lbd_{i}\lbd_{j}},\ &j\neq i,\\
\ \ 0,\ &j=i.
\end{cases}
$$
We denote that
$$\mathcal{Q}_{i}(\eta):=\frac{P_{m}\lbd_{i}}{\hat{f}}\Big(\tilde{A}_{i}+\tilde{B}_{i}+\tilde{C}_{i}+\tilde{D}_{i}-\left( 1+\frac{\eta}{m} \right)\tilde{E}_{i}\Big).$$
By a direct computation, we obtain
\begin{align}\label{k=n}
\mathcal{Q}_{i}(\eta)&=\lbd_{i}^{m}\left[(W\hat{f}-1)(\sum_{j}\frac{b_{jj,i}}{\lbd_{j}})^{2}+\sum_{j}(\frac{b_{jj,i}}{\lbd_{j}})^{2}\right]+(m-1)\sum_{j}\lbd_{j}^{m}(\frac{b_{jj,i}}{\lbd_{j}})^{2}\notag\\
&+2\sum_{j\neq i}(\lbd_{j}^{m}+O_{ij})(\frac{b_{jj,i}}{\lbd_{j}})^{2}-\frac{m+\eta}{P_{m}}(\sum_{j}\lbd_{j}^{m}\frac{b_{jj,i}}{\lbd_{j}})^{2},
\end{align}
where $O_{ij}:=\lbd_{i}\lbd_{j}\sum_{r=0}^{m-2}\lbd_{j}^{m-2-r}\lbd_{i}^{r}$ for $i\neq j$.

\vskip 2mm
\noindent
{\bf (1). Lemma \ref{h1} for the case $k=n$.} When $\eta=1$ and $i\neq 1$, let $m\geq 2$ and $W\hat{f}\geq 1$, and we get that
\begin{align*}
\mathcal{Q}_{i}(1)\geq \sum_{j}o_{j}(\frac{b_{jj,i}}{\lbd_{j}})^{2}-\frac{m+1}{P_{m}}(\sum_{j}\lbd_{j}^{m}\frac{b_{jj,i}}{\lbd_{j}})^{2},
\end{align*}
where
$$o_{i}:=m\lbd_{i}^{m},\quad o_{j}:=\lbd_{i}^{m}+(m+1)\lbd_{j}^{m}+2O_{ij}\ \text{for}\ j\neq i.$$
By applying the Cauchy-Schwarz inequality, we have
$$\Big(\sum_{j}\lbd_{j}^{m}\frac{b_{jj,i}}{\lbd_{j}}\Big)^{2}\leq \Big(\sum_{j}\frac{\lbd_{j}^{2m}}{o_{j}}\Big)\Big(\sum_{j}o_{j}(\frac{b_{jj,i}}{\lbd_{j}})^{2}\Big).$$
So we only need to prove
$$\sum_{j}\frac{\lbd_{j}^{2m}}{o_{j}}\leq \frac{P_{m}}{m+1}.$$
A direct computation gives,
\begin{align*}
 \frac{P_{m}}{m+1}-\sum_{j}\frac{\lbd_{j}^{2m}}{o_{j}}&=\frac{\lbd_{i}^{m}}{m+1}-\frac{\lbd_{i}^{m}}{m}+\sum_{j\neq i}(\frac{\lbd_{j}^{m}}{m+1}-\frac{\lbd_{j}^{2m}}{o_{j}})\\
 &=-\frac{\lbd_{i}^{m}}{m(m+1)}+\frac{1}{m+1}\sum_{j\neq i}\frac{\lbd_{j}^{m}(o_{j}-(m+1)\lbd_{j}^{m})}{o_{j}}\\
 &=\frac{1}{m+1}\Big(\sum_{j\neq i}\frac{\lbd_{j}^{m}(\lbd_{i}^{m}+2O_{ij})}{o_{j}}-\frac{\lbd_{i}^{m}}{m}\Big).
\end{align*}
Hence, for $j=1$, we obtain
\begin{align*}
m\lbd_{1}^{m}&(\lbd_{i}^{m}+2O_{i1})-\lbd_{i}^{m}\Big(\lbd_{i}^{m}+(m+1)\lbd_{1}^{m}+2O_{i1}\Big)\\
=&2O_{i1}(m\lbd_{1}^{m}-\lbd_{i}^{m})-\lbd_{i}^{m}(\lbd_{1}^{m}+\lbd_{i}^{m})\\
\geq& 2\lbd_{i}^{m}\lbd_{1}^{m}-\lbd_{i}^{m}(\lbd_{1}^{m}+\lbd_{i}^{m})\geq0,
\end{align*}
yields $\mathcal{Q}_{i}(1)\geq 0$.
\vskip 2mm
\noindent
{\bf (2). Lemma \ref{h2} for the case $k=n$.}
Next, we show that if $W\hat{f}\geq 2$, $0< \eta\leq \frac{1}{n}$, and $m\geq 2$, we have
$$\mathcal{Q}_{i}(\eta)\geq 0,\quad \text{for every}\quad i=1,\cdots,n.$$
Using the Cauchy-Schwarz inequality again,
\begin{align}\label{k=n--1}
\Big(\sum_{j}\lbd_{j}^{m}\frac{b_{jj,i}}{\lbd_{j}}\Big)^{2}\leq P_{m}\sum_{j}\lbd_{j}^{m}(\frac{b_{jj,i}}{\lbd_{j}})^{2}.
\end{align}
Combining with \eqref{k=n} and \eqref{k=n--1}, we have that
\begin{align*}
\mathcal{Q}_{i}(\eta)&\geq \lbd_{i}^{m}\left[(W\hat{f}-1)(\sum_{j}\frac{b_{jj,i}}{\lbd_{j}})^{2}+\sum_{j}(\frac{b_{jj,i}}{\lbd_{j}})^{2}\right]+(m-1)\sum_{j}\lbd_{j}^{m}(\frac{b_{jj,i}}{\lbd_{j}})^{2}\\
&\quad+2\sum_{j\neq i}(\lbd_{j}^{m}+O_{ij})(\frac{b_{jj,i}}{\lbd_{j}})^{2}-(m+\eta)(\sum_{j}\lbd_{j}^{m}\frac{b_{jj,i}}{\lbd_{j}})^{2}\\
&=\lbd_{i}^{m}\left[(W\hat{f}-1)(\sum_{j}\frac{b_{jj,i}}{\lbd_{j}})^{2}+\sum_{j\neq i}(\frac{b_{jj,i}}{\lbd_{j}})^{2}-\eta(\frac{b_{jj,i}}{\lbd_{i}})^{2}\right]\\
&\quad+(1-\eta)\sum_{j\neq i}\lbd_{j}^{m}(\frac{b_{jj,i}}{\lbd_{j}})^{2}+2\sum_{j\neq i}O_{ij}(\frac{b_{jj,i}}{\lbd_{j}})^{2}.
\end{align*}
Since $$(\sum_{j}\frac{b_{jj,i}}{\lbd_{j}})^{2}+\sum_{j\neq i}(\frac{b_{jj,i}}{\lbd_{j}})^{2}\geq \frac{(\frac{b_{jj,i}}{\lbd_{i}})^{2}}{n},$$
and therefore
\begin{align*}
(W\hat{f}-1)&(\sum_{j}\frac{b_{jj,i}}{\lbd_{j}})^{2}+\sum_{j\neq i}(\frac{b_{jj,i}}{\lbd_{j}})^{2}-\eta(\frac{b_{jj,i}}{\lbd_{i}})^{2}\\
&\geq (\sum_{j}\frac{b_{jj,i}}{\lbd_{j}})^{2}+\sum_{j\neq i}(\frac{b_{jj,i}}{\lbd_{j}})^{2}-\eta(\frac{b_{jj,i}}{\lbd_{i}})^{2}\\
&\geq (\frac{1}{n}-\eta)(\frac{b_{jj,i}}{\lbd_{i}})^{2}\geq 0.
\end{align*}
This completes the proof.
\end{proof}

\subsection{Proof of Theorem \ref{hxy66}.} We now complete the proof of Theorem \ref{hxy66}.

We first prove part $(1)$. Suppose that $p>q$ and $p\geq1$. Combining Lemmas \ref{lem3}, \ref{hxy1}, \ref{hxy2} and Corollary \ref{coro1} with Theorem \ref{thm2}, we obtain
$$\underset{\Ctheta}{\min}\ h\geq c,\quad\text{and}\quad ||h||_{C^{2}(\Ctheta)}\leq C.$$
The regularity theory for fully nonlinear second-order uniformly elliptic equations with oblique derivative boundary conditions \cite{LT}*{Theorem 1.1} then gives a $C^{2,\al}$ estimate for $h$, where $\al \in (0,1)$, together with the higher-order estimates in \eqref{hxy521}.

We next consider part $(2)$. If $1<q<p\leq q+1$, Eq. \eqref{eq:capillary-Lp-dual-CM} implies that the rescaled function
$$\tilde{h}:=\frac{h}{\underset{\Ctheta}{\min}\ h}$$
satisfies
\begin{equation*}
		\renewcommand{\arraystretch}{1.5}
		\left\{
		\begin{array}{rll}
			\sigma_k\big(\nabla^{2}\tilde{h}+\tilde{h}\sigma\big)&=(\underset{\Ctheta}{\min}\ h)^{p-q}f\tilde{h}^{p-1}\big(\tilde{h}^2+|\nabla \tilde{h}|^2\big)^{\frac{k+1-q}{2}},& \text{in }\Ctheta,\\
		\nabla_{\mu}\tilde{h}&=\cot\theta\, \tilde{h},& \text{on }\partial\Ctheta,
			\end{array}	
			\right.
	\end{equation*}

Suppose that $h$ attains its minimum and maximum at $\xi_{-}$ and $\xi_{+}$, respectively. Since $0<\theta<\pi/2$, the spherical cap $\Ctheta$ is geodesically convex. Let $\gamma$ be a unit-speed minimizing geodesic in $\Ctheta$ from $\xi_{-}$ to $\xi_{+}$. By Lemma \ref{hxy3},
\begin{align*}
\log\frac{\underset{\Ctheta}{\max}\ h}{\underset{\Ctheta}{\min}\ h}
&=\int_{0}^{\text{dist}(\xi_{-},\xi_{+})}\langle\nabla\log h(\gamma(s)),\dot\gamma(s)\rangle\,\mathrm{d}s\\
&\leq\operatorname{diam}(\Ctheta)~\underset{\Ctheta}{\max}|\nabla \log h|\leq C.
\end{align*}
This yields
$$1\leq \tilde{h}\leq\frac{\underset{\Ctheta}{\max}\ h}{\underset{\Ctheta}{\min}\ h}\leq C \quad \text{and} \quad |\nabla \tilde{h}|\leq \tilde{h}|\nabla \log h|\leq C.$$
Applying Theorem \ref{thm2} to $\tilde{h}$ yields
$$||\tilde{h}||_{C^{2}(\Ctheta)}\leq C.$$
The higher-order estimates \eqref{hxyhxy} now follow as in part $(1)$.
 \hfill{$\square$}

\section{Proof of main theorems}\label{sec 4}

When $1\leq k\leq n-1$, admissible solutions of Eq. \eqref{eq:capillary-Lp-dual-CM} are not necessarily strictly convex. The following theorem is to show that, under suitable assumptions on the function $f$, any convex solution of Eq. \eqref{eq:capillary-Lp-dual-CM} is in fact strictly convex. The proof combines the arguments in Hu and Cabezas-Moreno \cite{CabezasMorenoHu2025} with Hu, Ivaki and Scheuer \cite{HuIvakiScheuer2025}.
\begin{lemma}\label{full rank}
Let $1\leq k\leq n-1$, $p\geq 1$, $q\in\R$, and $\theta\in(0,\frac{\pi}{2})$. Suppose that $f$ is a positive smooth function satisfying conditions \eqref{23} and \eqref{24}. If $h$ is a positive smooth solution of Eq.\eqref{eq:capillary-Lp-dual-CM} such that $\nabla^2h+h\sigma\geq 0$ on $\Ctheta$, then $\nabla^2h+h\sigma$ is positive definite on $\Ctheta$.
\end{lemma}

\begin{proof}Define
$$F=\sigma_{k}^{\frac{1}{k}},\quad \tilde{f}=\hat{f}^{1/k}=(fh^{p-1}\rho^{k+1-q})^{1/k}.$$
Denote by $\lambda_{1}$ the smallest eigenvalue of $b[h]$. Suppose $\lambda_{1}=0$ somewhere in the interior of $\Ctheta$. By \cite{CabezasMorenoHu2025}*{Theorem 3.4, (3.20)} and condition \eqref{23}, there exists a positive constant $c$ depending on $k, p,q,||h||_{C^{3}(\Ctheta)}$ and $||f||_{C^{2}(\Ctheta)}$, such that
$$L[\lambda_{1}]:=F^{ij}\nabla^{2}_{ij}\lambda_{1}-c(\lambda_{1}+|\nabla \lambda_{1}|)\leq 0$$
in the viscosity sense. By the strong maximum principle, $\lambda_{1}\equiv0$ in $\Ctheta$. However, $\nabla_{\mu}h|_{\pa\Ctheta}=\cot\theta\, h>0$ implies that at the point where $h$ attains its minimum, we have $\lambda_{1}>0$. This is a contradiction.

Choose an orthonormal frame $\{e_{i}\}_{i=1}^{n}$ at $\xi_{0}$ such that
$$e_{1}=\mu,\quad e_{\al}\in T_{\xi_{0}}\pa \Ctheta.$$
We next assume that $\lambda_{1}(\xi_{0}) = 0$ for some $\xi_{0} \in \partial\Ctheta$ while $\lambda_{1} > 0$ in the interior. Then, by following the proof of \cite{HuIvakiScheuer2025}*{Theorem 3.1} with minor modifications, we obtain a contradiction. For completeness, we sketch the argument below.

\noindent
{\bf{Step 1}}. We prove
$$b_{ii}(\xi_{0})=0\quad\Longrightarrow\quad\nabla_{\mu}b_{ii}(\xi_{0})\geq 0.$$
By the Weingarten equation and Gauss formula, we have
\begin{equation}\label{GW}
b_{\al\al,\mu}=\cot\theta(b_{\mu\mu}-b_{\al\al}),
\end{equation}
where we used the fact $h_{\al \mu}=0$.
For $i=\al\geq2$, the claim follows immediately from \eqref{GW}. For $i=1$, recall that \eqref{eq:capillary-Lp-dual-CM} is equivalent to
\begin{equation}\label{feq}
F(b[h])=\tilde{f}\quad \text{in}\ \Ctheta.
\end{equation}
Differentiating \eqref{feq} in the $\mu$-direction gives
$$\sum_{i}F^{ii}b_{ii,\mu}=\tilde{f}_{\mu}.$$
By \eqref{GW}, we obtain
\begin{equation}\label{25}
F^{\mu\mu}b_{\mu\mu,\mu}=\tilde{f}_{\mu}+\sum_{\al}F^{\al\al}\cot\theta(b_{\al\al}-b_{\mu\mu}).
\end{equation}
Note that $b_{\mu\mu}(\xi_{0})=0$, and by the $1$-homogeneity of $F$, we get
$$\sum_{\al}F^{\al\al}b_{\al\al}=F=\tilde{f}.$$
Evaluating \eqref{25} at $\xi_{0}$ yields
$$b_{\mu\mu,\mu}=\frac{\tilde{f}_{\mu}+\tilde{f}\cot\theta}{F^{\mu\mu}}.$$
Thus,
$$
\frac{\tilde{f}_{\mu}+\tilde{f}\cot\theta}{F^{\mu\mu}}\geq0\Longleftrightarrow \nabla_{\mu}(\log \tilde{f})\geq -\cot\theta.
$$
Using properties (1) and (2) in Lemma \ref{lem1}, we obtain
\begin{align*}
&\frac{1}{k}\nabla_{\mu}\Big(\log f+(p-1)\log h+(k+1-q)\log \rho\Big)\\
&=\frac{1}{k}\Big(\frac{\nabla_{\mu}f}{f}+(p-1)\cot\theta+(k+1-q)\frac{\cot\theta hb_{\mu\mu}}{\rho^{2}}\Big)\\
&=\frac{1}{k}\Big(\frac{\nabla_{\mu}f}{f}+(p-1)\cot\theta\Big).
\end{align*}
Therefore
$$\nabla_{\mu}(\log \tilde{f})\geq -\cot\theta\Longleftrightarrow \nabla_{\mu}\log f\geq -(p+k-1)\cot\theta\quad \text{on}\ \pa\Ctheta,$$
which is precisely \eqref{24}.

\noindent
{\bf{Step 2}}. Consider an interior ball $B_{\tilde{R}}(x_0)\subset\Ctheta$ touching $\pa\Ctheta$ at $\xi_{0}$. Define an annular region $A_{\tilde{R},\tilde{\rho}}=B_{\tilde{R}}(x_0)\setminus  \operatorname{int}(B_{\tilde{\rho}}(x_0))$ for some $0<\tilde{\rho}<\tilde{R}$.
For $x\in \Ctheta$, let $\tilde{r}(x)=\operatorname{dist}(x,x_0)$ denote the distance from $x$ to $x_0$. We define
$$w(x)=e^{-\al \tilde{R}^2}-e^{-\al \tilde{r}(x)^2},$$
for a constant $\al>0$ to be determined later. For any $x\in A_{\tilde{R},\tilde{\rho}}$, the distance function $\tilde{r}(x)$ satisfies
$$
\nabla^2_{ij}\tilde{r}(x)=\cot \tilde{r}(x) (\sigma _{ij}-\nabla_{i}\tilde{r} \nabla_{j} \tilde{r}).
$$
Therefore,
\begin{align*}
\nabla_{i} w &=2\al \tilde{r} \nabla_{i}\tilde{ r} e^{-\al \tilde{r}^2},\\
\nabla^2_{ij} w&=-4\al^2 \tilde{r}^2 \nabla_{i} \tilde{r}  \nabla_{j} \tilde{r} e^{-\al \tilde{r}^2}+2\al \nabla_{i}\tilde{r} \nabla_{j} \tilde{r} e^{-\al \tilde{r}^2}+2\al \tilde{r} \nabla^2_{i,j} \tilde{r} e^{-\al \tilde{r}^2} \\
&=e^{-\al \tilde{r}^2}\Big( (-4\al^2 \tilde{r}^2+2\al-2\al \tilde{r} \cot\tilde{r}) \nabla_{i} \tilde{r} \nabla_{j} \tilde{r}+2\al \tilde{r} \cot \tilde{r} \sigma_{ij}\Big),
\end{align*}	
and
\begin{align*}
L[w]= &~e^{-\al \tilde{r}^2}(-4\al^2 \tilde{r}^2+2\al-2\al \tilde{r}\cot \tilde{r})|\nabla \tilde{r}|_{\dot{F}}^2\\
 &~+2\al \tilde{r}\cot \tilde{r} e^{-\al \tilde{r}^2}\operatorname{tr}(\dot{F})-2 c \al \tilde{r} e^{-\al \tilde{r}^2}+c(e^{-\al \tilde{r}^2}-e^{-\al \tilde{R}^2}).
\end{align*}

Assume that $\lambda \delta_{ij} \leq F^{ij} \leq \Lambda \delta_{ij}$ in $B_{\tilde{R}}(x_0)$, where $0<\lambda\leq \Lambda$ are some constants. Note that $0\leq \tilde{r}\cot \tilde{r} \leq 1$ for any $\tilde{\rho} \leq \tilde{r}\leq \tilde{R}$. Now take $\al>0$ sufficiently large so that in $A_{\tilde{R},\tilde{\rho}}$,
\begin{equation*}
L[w] \leq -e^{-\al \tilde{r}^2}\left[ (4\al^2\tilde{r}^2-2\al)\lambda -2\al n  \Lambda +2c\al \tilde{r}-c\right]-ce^{-\al \tilde{R}^2}<0.
\end{equation*}

\noindent
{\bf{Step 3}}. Since $\lambda_1(x)>0$ on $\partial B_{\tilde{\rho}}(x_0)$, there exists $\varepsilon>0$ such that
$$
\psi(x):=\lambda_1(x)+\varepsilon w(x) > 0
$$
on $\partial B_{\tilde{\rho}}(x_0)$. Note that $w = 0$ on $\partial B_{\tilde{R}}(x_0)$ and by our assumption there is no other point on $\partial B_{\tilde{R}}(x_0) \setminus \{\xi_{0}\}$ where $\lambda_1 = 0$; hence, $\psi > 0$ on $\partial B_{\tilde{R}}(x_0) \setminus \{\xi_{0}\}$. By the maximum principle for the viscosity supersolution $\psi$, we have $\psi\geq 0$ in the annulus, where $\psi$ is also the smallest eigenvalue of
$$
S_{ij}=b_{ij}+\varepsilon w \sigma_{ij}.
$$
Suppose that at $\xi_{0}$ the zero eigenvalue of $S(\xi_{0})$ is attained in direction $e_{i}$, i.e., we also have $b_{ii} = 0$.
Let $\ga$ be a unit-speed geodesic in the direction $-\mu$ and $e_i$ be parallel transported along $\ga$. Then from {\bf Step 1} we get
$$
0\leq \frac{d}{dt}\Big|_{t=0^+}S(\gamma(t))(e_i,e_i)=\nabla_{-\mu}b_{ii}+\varepsilon  \nabla_{-\mu}w\leq-2\alpha \varepsilon \tilde{R} <0,
$$
a contradiction.
\end{proof}

Using the a priori estimates, the full-rank lemma, and the continuity method, we establish existence and uniqueness of strictly convex solutions of Eq. \eqref{eq:capillary-Lp-dual-CM}.
\begin{proof}[\textbf{Proof of Theorem  \ref{thm1}}]
We first prove the existence and uniqueness assertion in part $(1)$. For a positive function $h$, set
$$v=\log h,\quad B[v]_{ij}=v_{ij}+v_{i}v_{j}+\sigma_{ij}.$$
The Robin condition becomes
$$\nabla_{\mu}v=\cot\theta\quad \text{on\ }\pa\Ctheta.$$
Using the homogeneity of $\sigma_{k}$, Eq. \eqref{eq:capillary-Lp-dual-CM} is equivalent to
\begin{equation*}
		\renewcommand{\arraystretch}{1.5}
		\left\{
		\begin{array}{rll}
		\Omega_{p,q,f}:&=\log \sigma_{k}(B[v])-\frac{k+1-q}{2}\log(1+|\nabla v|^{2})-(p-q)v-\log f=0,& \text{in }\Ctheta,\\
		\nabla_{\mu}v&=\cot\theta ,& \text{on }\partial\Ctheta.
			\end{array}	
			\right.
			\end{equation*}
The term
$$-(p-q)v$$
plays a key role. Since $p>q$, the zeroth-order coefficient of the linearization is strictly negative. Fix $p^*>k+1$ and define
$$p_{t}=(1-t)p^*+tp,\quad q_{t}=(1-t)(k+1)+tq,\quad 0\leq t\leq 1.$$
Then
$$p_{t}-q_{t}=(1-t)(p^*-k-1)+t(p-q)\geq \delta_{0}>0,$$
where $\delta_{0}:=\min\{p^*-k-1,p-q\}$. Now put
$$m_{t}:=p_{t}+k-1,\quad g_{t}:=(1-t)g_{0}+tg_{1},$$
where $g_{0}=\binom{n}{k}^{-\frac{1}{p^*+k-1}}\ell^{\tilde{\lambda}}$ with $\tilde{\lambda}=\frac{p^*-1}{p^*+k-1}$ and $g_{1}:=f^{-\frac{1}{p+k-1}}$. We first check that $g_{t}$ satisfies
$$\nabla^{2}g_{t}+g_{t}\sigma\geq 0\quad \text{in}\ \Ctheta,\quad \nabla_{\mu}g_{t}\leq \cot\theta g_{t}\quad\text{on}\ \pa\Ctheta.$$

A direct computation shows that
\begin{align*}
\nabla^{2}(\ell^{\tilde{\lambda}})+\ell^{\tilde{\lambda}}\sigma&=\Big(\tilde{\lambda}\ell^{\tilde{\lambda}-1}+(1-\tilde{\lambda})\ell^{\tilde{\lambda}}\Big)\sigma+\tilde{\lambda}(\tilde{\lambda}-1)\ell^{\tilde{\lambda}-2}\nabla \ell\otimes\nabla\ell\\
&\geq \tilde{\lambda}\ell^{\tilde{\lambda}-2}\Big(\ell \sigma-(1-\tilde{\lambda})\nabla \ell\otimes\nabla\ell\Big).
\end{align*}
Moreover,
$$\ell^{2}+|\nabla \ell|^{2}=2\ell-\sin^{2}\theta,\quad 0< \ell\leq \sin^{2}\theta,$$
and hence
$$\ell-|\nabla \ell|^{2}=\sin^{2}\theta-\ell+\ell^{2}.$$
It follows that
$$\nabla^{2}g_{0}+g_{0}\sigma>0.$$
On $\pa\Ctheta$,
$$\nabla_{\mu}g_{0}=\tilde{\lambda}\cot\theta g_{0}<\cot\theta g_{0},$$
and both conditions are preserved by $g_t$ for every $t\in[0,1]$. This allows us to apply the full-rank theorem at a limit point. We consider the following family of equations
\begin{equation}
		\renewcommand{\arraystretch}{1.5}
		\left\{
		\begin{array}{rll}\label{520}
			\sigma_k\big(\nabla^{2}h+h\sigma\big)&=f_{t}h^{p_{t}-1}\big(h^2+|\nabla h|^2\big)^{\frac{k+1-q_{t}}{2}},& \text{in }\Ctheta,\\
		\nabla_{\mu}h&=\cot\theta\, h,& \text{on }\partial\Ctheta,
			\end{array}	
			\right.
	\end{equation}
where $f_{t}:=g_{t}^{-m_{t}}$. Fix $0<\al<1$ and define
$$\mathcal{V}_{\al}:=\{v\in C^{2,\al}(\Ctheta)\ :\ B[v]>0\quad\text{and}\quad \nabla_{\mu}v=\cot\theta\quad\text{on}\ \pa\Ctheta\}.$$
For $t\in [0,1]$, let
$$\Omega_{t}[v]:=\log \sigma_{k}(B[v])-\frac{k+1-q_{t}}{2}\log(1+|\nabla v|^{2})-(p_{t}-q_{t})v-\log f_{t}$$
and define
$$\mathcal{G}(t,v)=(\Omega_{t}[v],\nabla_{\mu}v-\cot\theta).$$
Let
$$\mathcal{I}:=\{t\in [0,1]\ :\ \exists v\in \mathcal{V}_{\al}\quad\text{with}\quad \mathcal{G}(t,v)=(0,0)\}.$$
\noindent
{\bf Step 1. Non-emptiness.}

At $t=0$, $h_{0}=\ell$, equivalently $v_{0}=\log \ell$ solves \eqref{520}. Hence $0\in \mathcal{I}$.

\noindent
{\bf Step 2. Openness.}

Let $v\in \mathcal{V}_{\al}$. Define $G(B[v]):=\log \sigma_{k}(B[v])$, and
$$G^{ij}:=\frac{\pa \log \sigma_{k}}{\pa B_{ij}}(B[v])=\frac{\sigma_{k}^{ij}(B[v])}{\sigma_{k}(B[v])}.$$
Then the linearization in a variation $\psi$ is
\begin{equation*}
D_{v}\Omega_{t}[v](\psi)=G^{ij}\psi_{ij}+\Big(2G^{ij}v_{i}-\frac{k+1-q_{t}}{1+|\nabla v|^{2}}v^{j}\Big)\psi_{j}-(p_{t}-q_{t})\psi.
\end{equation*}
The linearized boundary condition is
$$\nabla_{\mu}\psi=0,\quad \text{on}\ \pa\Ctheta.$$
Since $B[v]>0$ on the compact domain, this operator is uniformly elliptic, and the zeroth-order coefficient $p_{t}-q_{t}$ is strictly negative. Hence the maximum principle and the Hopf boundary lemma give a trivial kernel for the homogeneous oblique problem. Standard Schauder theory for uniformly elliptic equations with uniformly oblique boundary condition then gives that
$$\psi\mapsto(D_{v}\Omega_{t}[v](\psi),\nabla_{\mu}\psi)$$
is an isomorphism
$$C^{2,\al}(\Ctheta)\mapsto C^{0,\al}(\Ctheta)\times C^{1,\al}(\pa\Ctheta).$$
The Banach-space implicit function theorem therefore shows that $\mathcal{I}$ is open.

\noindent
{\bf Step 3. Closedness.}

Let $t_{j}\rightarrow t_{\infty}$, and let $h_{j}$ be the corresponding positive, smooth, strictly convex solutions. The family $(p_{t},q_{t},f_{t})$ is compact in the parameter--data space, and the gap $p_{t}-q_{t}$ has a uniform positive lower bound. The a priori and higher-order estimates for $p>q$ therefore yield, after passing to a subsequence,
$$h_{j}\rightarrow h_{\infty}\quad \text{in}\ C^{4,\ga}(\Ctheta)$$
for some $\ga\in (\al,1)$, with $h_{\infty}>0$. Since $b[h_{j}]>0$, we have $b[h_{\infty}]\geq 0$. Passing to the limit in \eqref{520} yields
$$\sigma_{k}(b[h_{\infty}])=f_{t_{\infty}}h_{\infty}^{p_{t_{\infty}}-1}(h_{\infty}^{2}+|\nabla h_{\infty}|^{2})^{\frac{k+1-q_{t_{\infty}}}{2}}>0.$$
The function $f_{t_{\infty}}$ satisfies the hypotheses of the full-rank theorem because $g_{t_{\infty}}=f_{t_{\infty}}^{-1/({p_{t_{\infty}}+k-1})}$ satisfies \eqref{23} and \eqref{24}. The Lemma \ref{full rank} therefore gives $b[h_{\infty}]>0$.
Hence $h_{\infty}$ is a positive, smooth, strictly convex solution at $t_{\infty}$, so $t_{\infty}\in \mathcal{I}$. Thus $\mathcal{I}$ is closed. Setting $t=1$ gives a smooth, positive, strictly convex solution of Eq. \eqref{eq:capillary-Lp-dual-CM}.

It remains to prove uniqueness. Suppose that $h_{1}$ and $h_{2}$ are two positive strictly convex solutions of Eq. \eqref{eq:capillary-Lp-dual-CM}. Set $\mathcal{K}:=\frac{h_{1}}{h_{2}}$; then $\nabla_{\mu}\mathcal{K}=0$ on $\pa\Ctheta$. Suppose that $\mathcal{K}$ attains its maximum at $\xi_{0}$. At this point,
$$0=\nabla\log \mathcal{K}=\frac{\nabla h_{1}}{h_{1}}-\frac{\nabla h_{2}}{h_{2}},$$
and
$$0\geq\nabla^{2}\log\mathcal{K}=\frac{\nabla^{2}h_{1}}{h_{1}}-\frac{\nabla^{2}h_{2}}{h_{2}},$$
which imply
$$\sigma_{k}\Big(\frac{\nabla^{2}h_{1}+h_{1}\sigma}{h_{1}}\Big)\leq \sigma_{k}\Big(\frac{\nabla^{2}h_{2}+h_{2}\sigma}{h_{2}}\Big).$$
Therefore,
\begin{align*}
1=\frac{f(\xi_{0})}{f(\xi_{0})}=\frac{(h_{1}^{2}+|\nabla h_{1}|^{2})^{-\frac{k+1-q}{2}}h_{1}^{1-p}\sigma_{k}(\nabla^{2}h_{1}+h_{1}\sigma)}{(h_{2}^{2}+|\nabla h_{2}|^{2})^{-\frac{k+1-q}{2}}h_{2}^{1-p}\sigma_{k}(\nabla^{2}h_{2}+h_{2}\sigma)}\leq \mathcal{K}(\xi_{0})^{q-p}.
\end{align*}
Since $p>q$, it follows that $\underset{\Ctheta}{\max}\ \mathcal{K}=\mathcal{K}(\xi_{0})\leq 1$. An analogous argument at a minimum gives $\underset{\Ctheta}{\min}\ \mathcal{K}\geq 1$. Consequently, $h_{1}\equiv h_{2}$.

We now prove part $(2)$, where $p=q>1$. For $\ep\in(0,1)$, set
$$f_{\ep}:=f^{\frac{p+\ep+k-1}{p+k-1}}.$$
Then $f_{\ep}\to f$ smoothly, and
$$f_{\ep}^{-\frac{1}{p+\ep+k-1}}=f^{-\frac{1}{p+k-1}}.$$
Thus the convexity and boundary conditions in the theorem hold for the approximating data. Consider the approximating equation
\begin{equation}\label{lml1}
		\renewcommand{\arraystretch}{1.5}
		\left\{
		\begin{array}{rll}
			\sigma_k\big(\nabla^{2}h+h\sigma\big)&=f_{\ep}h^{p-1+\ep}\big(h^2+|\nabla h|^2\big)^{\frac{k+1-p}{2}},& \text{in }\Ctheta,\\
		\nabla_{\mu}h&=\cot\theta\, h,& \text{on }\partial\Ctheta.
			\end{array}	
			\right.
	\end{equation}
The preceding argument provides a positive solution $h_{\ep} \in C^{4,\al}(\Ctheta)$ of Eq. \eqref{lml1}. Set
$$m_{\ep}:=\min_{\Ctheta}h_{\ep},\qquad
\tau_{\ep}:=m_{\ep}^{\ep},\qquad
\tilde{h}_{\ep}:=\frac{h_{\ep}}{m_{\ep}}.$$
Write $u_{\ep}=h_{\ep}/\ell$. The maximum and minimum comparisons in the proof of Lemma \ref{lem3}, applied with $(p+\ep,p,f_{\ep})$, give
$$\bigl(\max_{\Ctheta}u_{\ep}\bigr)^{\ep}\leq C,
\qquad \bigl(\min_{\Ctheta}u_{\ep}\bigr)^{\ep}\geq c.$$
The constants $c,C>0$ are independent of $\ep$. Since $\ell$ has a fixed positive upper and lower bound, it follows that $c\leq\tau_{\ep}\leq C$.
It follows from \eqref{lml1} that $\tilde{h}_{\ep}$ satisfies
\begin{equation*}
		\renewcommand{\arraystretch}{1.5}
		\left\{
		\begin{array}{rll}
			\sigma_k\big(\nabla^{2}\tilde h_\ep+\tilde h_\ep\sigma\big)&=\tau_\ep f_\ep\tilde h_\ep^{p-1+\ep}\big(\tilde h_\ep^{2}+|\nabla\tilde h_\ep|^2\big)^{\frac{k+1-p}{2}},& \text{in }\Ctheta,\\
		\nabla_{\mu}\tilde h_\ep&=\cot\theta\,\tilde h_\ep,& \text{on }\partial\Ctheta.
			\end{array}	
			\right.
	\end{equation*}
The uniform bounds on $\tau_\ep$ and $f_\ep$, together with Theorem \ref{hxy66}, give
$$||\tilde{h}_{\ep}||_{C^{4,\al}(\Ctheta)}\leq C,$$
where $C$ is independent of $\ep$. Passing to a subsequence, we have $\tilde h_{\ep_i}\to h$ in $C^{4,\al_0}(\Ctheta)$ for $0<\al_0<\al$ and $\tau_{\ep_i}\to\tau>0$. Thus $(h,\tau)$ solves
\begin{equation*}
		\renewcommand{\arraystretch}{1.5}
		\left\{
		\begin{array}{rll}
			\sigma_k\big(\nabla^{2}h+h\sigma\big)&=\tau fh^{p-1}\big(h^2+|\nabla h|^2\big)^{\frac{k+1-p}{2}},& \text{in }\Ctheta,\\
		\nabla_{\mu}h&=\cot\theta\, h,& \text{on }\partial\Ctheta.
			\end{array}	
			\right.
	\end{equation*}
The limit satisfies $b[h]\geq0$ and $\sigma_k(b[h])>0$. The Lemma \ref{full rank} gives $b[h]>0$. This proves existence in part $(2)$. We next establish uniqueness when $p=q$.

\noindent
{\bf Step 1. Uniqueness of the constant $\tau$.} Assume $(h_{1},\tau_{1})$ and $(h_{2},\tau_{2})$ are two solution pairs. They satisfy
\begin{equation*}
		\renewcommand{\arraystretch}{1.5}
		\left\{
		\begin{array}{rll}
			\log \sigma_{k}(B[v])-\frac{k+1-p}{2}\log(1+|\nabla v|^{2})&=\log \tau+\log f,& \text{in }\Ctheta,\\
		\nabla_{\mu}v&=\cot\theta ,& \text{on }\partial\Ctheta.
			\end{array}	
			\right.
	\end{equation*}
At a maximum point of $v_{1}-v_{2}$ on the closed cap, the boundary-extremum argument used in the uniqueness proof for part $(1)$ gives $\tau_{1}\leq \tau_{2}$. Applying the analogous argument at a minimum yields the reverse inequality; hence $\tau_{1}=\tau_{2}$.

\noindent
{\bf{Step 2}. Uniqueness of the solution.} Set $\tilde{w}:=v_{1}-v_{2}$ and define
$$B_{t}=(1-t)B[v_{2}]+tB[v_{1}],\quad a^{ij}=\int_{0}^{1}G^{ij}(B_{t})dt,\quad u_{t}=(1-t)\nabla v_{2}+t\nabla v_{1}.$$
Subtracting the two equations and applying the fundamental theorem of calculus to $G(B_{t})$ and the gradient term gives
\begin{equation*}
\renewcommand{\arraystretch}{1.5}
\left\{
\begin{array}{rll}
a^{ij}\tilde{w}_{ij}+b^{i}\tilde{w}_{i}&=0,\quad& \text{in }\Ctheta,\\
\tilde{w}_{\mu}&=0,\quad&\text{on }\pa\Ctheta,
\end{array}	
\right.
\end{equation*}
where
$$b^{i}=a^{ij}\Big((v_{1})_{j}+(v_{2})_{j}\Big)-(k+1-p)\Big(\int_{0}^{1}\frac{(u_{t})_{i}}{1+|u_{t}|^{2}}dt\Big).$$
These coefficients are bounded, and $a^{ij}$ is uniformly positive definite on the compact cap. The strong maximum principle and the boundary point lemma therefore imply that $\tilde{w}$ is constant. Thus $h_{1}=ch_{2}$ for some $c>0$. If both solutions satisfy $\underset{\Ctheta}{\min}\ h=1$, then $h_{1}\equiv h_{2}$. This completes the proof.
 \end{proof}

 \begin{proof}[\textbf{Proof of Theorem  \ref{thm11}}]
Theorem \ref{thm11} then follows directly from Theorem \ref{thm1}, which establishes the strict convexity of solutions to Eq. \eqref{eq:capillary-Lp-dual-CM} without imposing either the symmetry condition \eqref{even} or
the convexity conditions \eqref{23} and \eqref{24} on $f$.

 \end{proof}
\section*{Conflict of Interest Statement}
The authors declare that they have no conflict of interest.

\section*{Data Availability Statement}
No datasets were generated or analyzed during the current study.

\section*{References}


\end{document}